\documentclass{article}
\pdfoutput=1
\usepackage{amsmath}
\usepackage{graphicx}
\newtheorem{remark}{Remark}

\begin{document}

\title{Higher order semi-implicit schemes for linear advection equation on Cartesian grids with numerical stability analysis}

\author{Peter Frolkovi\v{c}, Karol Mikula\thanks{This work was supported by grants VEGA 1/0728/15, VEGA 1/0608/15 and APVV-15-0522.}\\
Slovak University of Technology in Bratislava, Faculty of Civil Engineering\\
Department of Mathematics and Descriptive Geometry\\
Radlinsk\'eho 11, 810 05 Bratislava, Slovak Republic\\
\tt{peter.frolkovic@stuba.sk, karol.mikula@stuba.sk} 
}

\maketitle

\begin{abstract}
A new class of semi-implicit numerical schemes for linear advection equation on Cartesian grids is derived that is inspired by so-called $\kappa$-schemes \cite{lee77,wes01} used with fully explicit discretizations for this type of problems. Opposite to fully explicit $\kappa$-scheme the semi-implicit variant is unconditionally stable in one-dimensional case and it preserves second order accuracy for dimension by dimension extension in higher dimensional cases. We discuss von Neumann stability conditions numerically for all numerical schemes. Using so-called Corner Transport Upwind extension \cite{col90,lev02} of two-dimensional semi-implicit scheme with a special choice of $\kappa$ parameters, a second order accurate method is obtained for which numerical unconditional stability can be shown for variable velocity and the third order accuracy can be proved for constant velocity. Several numerical experiments illustrate the properties of semi-implicit schemes for chosen examples.
\end{abstract}


\section{Introduction}

In this work we derive a new class of semi-implicit schemes for the numerical solution of model advection equation
\begin{equation}
\nonumber
\partial_t u(x,t) + \vec{V}(x) \cdot \nabla u(x,t) = 0 \,, \quad u(x,0) = u^0(x) \,.
\end{equation}
This type of equation is a part of many mathematical models that are used in applications like flow and transport problems \cite{pat80,bea87,wes01,rnt16} and level set methods \cite{set99,osh02} for tracking of interfaces \cite{sus94,fed02,xu06,gro11,fro12,flw16}.

The numerical schemes will be derived for Cartesian grids that are used especially if input geometry has already a regular structure \cite{set99,fb09,nat16}, or if geometric boundaries are treated implicitly like in Cartesian cut cell methods \cite{lev88b,ing03}, immersed interface methods \cite{lev94,li06,xu06}, ghost fluid methods \cite{fed99,fed00} and similar. 

Standard methods for numerical solutions of advection equation are so-called fully explicit numerical schemes \cite{lev02,osh03}. The main advantage is their simplicity as the numerical solution, once the scheme is constructed, can be obtained directly. 
On the other hand the well-known restriction of fully explicit schemes is the CFL stability condition on the choice of time step that depends, among other, on a length of grid step size. Although such restriction is not considered as a disadvantage in general, it appears critical if geometric boundaries are resolved only implicitly \cite{lev88b,lev94,fed99}. The presence of arbitrary small cut cells in such treatment can give locally arbitrary small grid size that results in unrealistic CFL restriction if no modification of numerical scheme is provided.

Even if a stability restriction on the size of time step is not critical due to space discretization, it can occur if the model to be solved is not of pure advection form, but it contains a PDE or PDEs with several coupled terms.
For such problems  the fully explicit discretization of advection decouples this term from others using e.g. operator splitting approach \cite{mar90,hun13} that can cause an inaccuracy or instability in numerical solutions \cite{val92,lev02,rnt16}. 

For these reasons, numerical schemes based on fully implicit time discretization are often considered in such situations \cite{kna05,hun13,rnt16}. Opposite to fully explicit variants, the fully implicit schemes can couple the unknown values of numerical solution in all discretized terms of PDEs, and avoid the time splitting error. The price to pay is that a system of algebraic equations (in general nonlinear) has to be solved to obtain the numerical solution.

The implicit schemes are supposed to have less restrictive stability requirements concerning the choice of time steps. A well-known representative example is the first order accurate upwind scheme for linear advection \cite{pat80,lev02} that is unconditionally stable concerning the choice of time step. 
Nevertheless, the scheme is producing very inaccurate numerical solutions in too many situations. Clearly, numerical schemes with improved accuracy and unconditional stability are of general interest for computational purposes.

Recently, several publications \cite{mo11,mou14,fmu14,fro16} are dealing with semi-implicit finite volume schemes that combine the both types of time discretization even for pure advection equation. The main idea is that the 
implicit time discretization is used only for numerical solution at inflow boundaries of computational cells.

The aim of this paper is to extend this approach in several aspects for Cartesian grids. Firstly we aim to derive such semi-implicit schemes in finite difference form. When compared with finite volume formulation this form is computationally less demanding and can be characterized directly concerning stability conditions. 
The main novel tool in the derivation is so-called partial Cauchy-Kowalevski (or Lax-Wendroff) procedure. In its standard form this procedure replaces the time derivatives of solution in Taylor series by using only space derivatives of solution obtained from the advection equation \cite{tor09}. In our approach we apply less steps of such procedure by allowing also mixed time-space derivatives.

Next we derive a class of semi-implicit schemes following the approach of fully explicit $\kappa$-scheme \cite{lee77,lee85,wes01} that includes several popular numerical schemes like Lax-Wendroff or Fromm scheme \cite{lev02}. Such parametric formulation of the scheme can be used to improve its accuracy in special cases or to optimize it using so-called limiters \cite{lev02,hir07,tor09}. 

To compare clearly the semi-implicit $\kappa$-scheme with an analogous fully implicit scheme, we derive the both variants in 1D.
The schemes are second order accurate for one-dimensional advection equation with variable velocity. Moreover, they offer a special (velocity dependent) value of $\kappa$ for which the scheme is third order accurate if the velocity is constant. Nevertheless, we can show that the semi-implicit variant has two important advantages when compared to the fully implicit variant.

Firstly, it preserves the second order accuracy when applied in several dimensional case using simple dimension by dimension extension of 1D scheme.  Such property is not valid for fully explicit $\kappa$-schemes \cite{lev02} and, analogously,  also not for fully implicit $\kappa$-scheme.

Secondly, the semi-implicit $\kappa$-scheme has better stability properties than the fully implicit variant. Note that when deriving stability conditions of all schemes in this paper we apply von Neumann stability analysis \cite{tre96,wes01,hir07} realized in numerical way as suggested in \cite{bil97a,bil97b}. Although such approach gives enough confidence in obtained stability results, such analysis can not be considered rigorous therefore we use the notion of numerical stability analysis or numerical stability conditions. 

For the one-dimensional form of semi-implicit $\kappa$-scheme we can report unconditional numerical stability for all interesting values of $\kappa$. This is not the case of fully implicit $\kappa$-scheme when such stability result is obtained only for a subset of values for $\kappa$. 

Concerning the numerical stability analysis of higher dimensional case we restrict in this paper only to two-dimensional case of semi-implicit $\kappa$-scheme that is already quite demanding. We found that such scheme has conditional numerical stability conditions, nevertheless they are significantly less restrictive than in the case of fully explicit schemes.

To further improve the accuracy and stability of two-dimensional semi-implicit $\kappa$-scheme we apply the idea of Corner Transport Upwind (CTU) scheme \cite{col90,lev02} by adding additional discretization terms to the scheme. The semi-implicit $\kappa$-scheme with such CTU extension using the special (velocity dependent) value of $\kappa$ has unconditional numerical stability. Moreover in the case of constant velocity vector this scheme is third order accurate.

The paper is organized as follows.  In section \ref{sec-1d} we begin with one-dimensional case where the fully explicit, fully implicit and semi-implicit $\kappa$-schemes are derived. In 
Section \ref{sec-2d} we discuss the properties of semi-implicit $\kappa$-scheme in higher dimensional case when obtained by simple dimension by dimension extension. Moreover the Corner Transport Upwind extension of such scheme in two-dimensional case is given. In section \ref{sec-num} several numerical experiments are presented that confirm the properties of semi-implicit schemes for chosen examples. Finally we conclude the results in section \ref{sec-con}.


\section{One dimensional case}
\label{sec-1d}

We begin with the derivation of numerical schemes for one dimensional advection equation written as
\begin{equation}
\partial_t u(x,t) + V(x) \partial_x u(x,t) = 0 \,, \quad u(x,0) = u^0(x) \,, \quad x \in R \,, \quad t \ge 0 \, . 
\label{adveq1d}
\end{equation}
Let $x_i$ be points of an uniform grid such that $x_i - x_{i-1} \equiv h$ with $h$ being a constant grid size.
Furthermore let $\tau>0$ be a given time step and $t^n=n \tau$, $n=0,1,\ldots$ . We use standard indexing for discrete values like $V_i=V(x_i)$ and so on. 

Our aim is to find approximate values $U_i^n$ such that $U_i^n \approx u_i^n$ where $u_i^n := u(x_i,t^n)$ 
and $U_i^0=u^0_i=u^0(x_i)$. 
All numerical schemes in this section will be based formally on a 
particular choice of coefficients in the following general scheme,
\begin{equation}
\label{all1d}
U_i^{n+1} + \sum_{k=-2}^{2} \alpha_{i k} U_{i+k}^{n+1} = U_i^n + \sum_{k=-2}^{2} \beta_{i k} U_{i+k}^{n} \,.
\end{equation}
The fully explicit form of (\ref{all1d}) will be given by $\alpha_{i k} \equiv 0$, analogously $\beta_{ i k} \equiv 0$ in the case of fully implicit form.
Our aim is to derive semi-implicit schemes that will have in general three consecutive nonzero values of coefficients $\alpha_{i k}$ and $\beta_{i k}$ in (\ref{all1d}). All presented numerical schemes will be characterized by an order of their accuracy and by conditions for their numerically indicated stability.

To check the order of accuracy for any particular scheme of (\ref{all1d}) we simply consider its truncation error that is obtained by replacing all numerical values $U_{i+k}^{n+1}$ and $U_{i+k}^{n}$  in (\ref{all1d}) with the exact values $u_{i+k}^{n+1}$ and $u_{i+k}^{n}$ that themselves are then expressed with Taylor series, see some standard textbooks on numerical analysis, e.g. \cite{hir07}.

To derive a stability condition of particular numerical scheme for (\ref{all1d}), we use the approach of von Neumann stability analysis, see e.g. \cite{tre96,wes01,hir07}. To do so one introduces a grid function $\epsilon_i^n=\epsilon(x_i,t^n)$ defined by  
\begin{equation}
\label{eps1d}
\epsilon(x,t)=\exp(-\lambda t) \exp(\imath x) \,, \quad x \in R \,, \,\, t \ge 0 \,,
\end{equation}
where $\imath$ is the imaginary number, and the parameter $\lambda$ shall be found. 
The values $\epsilon_i^n$ are supposed to fulfill the numerical scheme. Using trivial relations
\begin{equation}
\label{trivial}
\epsilon_{i\pm k}^n =\exp(\pm \imath k h) \epsilon_i^n \,, \quad
\epsilon_i^{n+1}=S \epsilon_i^n \,, \,\, S:=\exp(-\lambda \tau) \,, 
\end{equation}
where $S$ denotes the so-called amplification factor, the {von Neumann stability condition of numerical scheme} is derived by searching for conditions under which one has $|S| \le 1$ for all $h \in (-\pi,\pi)$ . Using (\ref{trivial}) in (\ref{all1d}) one obtains
\begin{equation}
\label{allstab}
S = \left(1 + \sum_{k=-2}^{2} \beta_{i k} \exp(\imath k h) \right) \left(1+ \sum_{k=-2}^{2} \alpha_{i k}  \exp(\imath k h) \right)^{-1} \,.
\end{equation}

Although the stability conditions for $|S| \le 1$ from (\ref{allstab}) can be derived using analytical methods  for some schemes \cite{tre96,wes01,hir07}, we apply an approach proposed and used in \cite{bil97a,bil97b} where such condition is found numerically. A straightforward approach is to compute the values $|S|$  for very large number of discrete values of $h$ and input parameters for particular numerical scheme \cite{bil97a,bil97b}. Additionally, we apply numerical optimization algorithms available in Mathematica \cite{mat16} to search for local maxima of $|S|$ for these initial guesses. The advantage of numerical stability analysis is that it can be applied to all numerical schemes studied in this paper including nontrivial two-dimensional variants later. We adopt the notion {numerical stability condition} using this approach as it can not be considered  rigorous especially when unconditional stability is numerically indicated for a particular method.

We present now the derivation of numerical schemes for fully explicit variants of (\ref{all1d})  in such way that it will be rather straightforward to extend it later for fully implicit and semi-implicit variants. In what follows we use shorter notations for the exacts values of derivatives by $\partial_t u_i^n := (\partial_t u)_i^n$ and so on. Analogous notation with capital letter $U$ is reserved for numerical approximations of derivatives.

An important role in our derivation will play the following parametric class of approximations $\partial_x^{\kappa} U_i^n \approx \partial_x u_i^n$ (and analogously later for $\partial_x^{\kappa} U_i^{n+1} \approx \partial_x u_i^{n+1}$)
\begin{eqnarray}
\label{kappagrad}
2 \, \partial_x^{\kappa} U_{i}^n := (1-\kappa) \, \partial_x^- U_i^n + (1+\kappa) \, \partial_x^+ U_i^n \,, 
\end{eqnarray}
where
$$
h \, \partial_x^- U_i^n := U_i^n-U_{i-1}^n \,, \quad h \, \partial_x^+ U_i^n := U_{i+1}^n-U_i^n \,,
$$
and the parameter $\kappa$ in (\ref{kappagrad})  is free to choose. A natural choice $\kappa \in [-1,1]$ gives a convex combination of standard one-sided finite difference approximations.

To derive the fully explicit schemes we express the exact solution using Taylor series in a forward manner for some integer $p$,
\begin{equation}
\label{ext}
u_i^{n+1} = u_i^n + \sum_{m=1}^p \frac{1}{m!} \tau^m \partial^m_t u_i^n  + \mathcal{O}(\tau^{p+1}) \,.
\end{equation}
Next the so-called Cauchy-Kowalewski procedure, also known as Lax-Wendroff procedure \cite{tor09}, is applied to (\ref{ext}) by replacing all time derivatives of $u$ by space derivatives of $u$ using the equation (\ref{adveq1d}). We can write it in the form
\begin{eqnarray}
\label{ck01}
\partial_t u_i^n = - V_i \partial_x u_i^n \,,\\
\label{ck02}
\partial_{tx} u_i^n = - \partial_x (V \partial_{x} u)_i^n \,,\\
\label{ck11}
\partial_{tt} u_i^n = - V_i \partial_{tx} u_i^n = V_i \partial_x (V \partial_{x} u)_i^n \,,
\end{eqnarray}
and use it in (\ref{ext}) for $p=2$ to obtain
\begin{eqnarray}
\label{accur1}
u_i^{n+1} =  u_i^n - \tau V_i \partial_x u_i^n + 0.5 \tau^2 V_i \partial_{x} (V \partial_x u)_i^n
+ \mathcal{O}(\tau^3) \,.
\end{eqnarray}

The fully explicit $\kappa$-scheme to solve (\ref{adveq1d}) is now obtained by applying proper (upwinded) finite difference approximations in (\ref{accur1}). To reach a truncation error corresponding to $2^{nd}$ order accurate schemes, the term multiplied by $\tau$ in (\ref{accur1}) must be approximated by $2^{nd}$ order accurate approximation, while for the term multiplied by $0.5 \tau^2$ a first order accurate approximation is sufficient.

Having this in mind we apply in (\ref{accur1}) the upwinded approximations
\begin{eqnarray}
\label{upwgrad1d}
V_i \partial_x u_i^n \approx
[V_i]^+\partial_x^- (U_i^n + 0.5 h \partial^{\kappa}_x U_i^n) +[V_i]^- \partial_x^+ (U_i^n - 0.5 h \partial^{\kappa}_x U_i^n)  \,, \\[1ex]\label{dxx}
V_i \partial_{x} (V \partial_x u)_i^n \approx \left( [V_i]^+ \partial^-_x + [V_i]^- \partial^+_x\right) (V_i \partial_x^{\kappa} U_i^n)
\end{eqnarray}
where $[V]^+=\max\{0,V\}$ and $[V]^-=\min\{0,V\}$. 
The fully explicit $\kappa$-scheme takes then the form
\begin{eqnarray}
\label{fexpl}
U_i^{n+1} = U_i^{n} -  \tau [V_i]^+ \partial_x^- \left( U_i^{n} + 0.5 (h-  \tau  V_i) \partial_x^{\kappa} U_i^n)\right) - \\[1ex] 
\nonumber
\tau [V_i]^- \partial_x^+ \left( U_i^{n} - 0.5 (h + \tau  V_i) \partial_x^{\kappa} U_i^n)\right) \,.
\end{eqnarray}

The scheme (\ref{fexpl})  gives in the case of constant velocity $V_i \equiv V$ the well-known particular variants \cite{wes01,lev02,tor09}, namely Lax-Wendroff for $\kappa = \hbox{sign}(V)$, Beam-Warming for $\kappa = - \hbox{sign}(V)$, and Fromm scheme for $\kappa = 0$. For variable velocity case the scheme (\ref{fexpl}) is used in \cite{fm07,fw09} in a finite volume context.

We present two particular variants of (\ref{fexpl}) for the case of variable velocity when $V_i>0$ to present them in a more clear way.
To do so we denote a nondimensional (in general signed) Courant numbers
\begin{equation}
\label{courant}
\mathcal{C}_i = \tau V_i / h \,,
\end{equation}
then for $\kappa = 1$ the scheme (\ref{fexpl}) takes the form
\begin{eqnarray}
\label{fexplLW}
U_i^{n+1} = U_i^{n} -  \mathcal{C}_i \left( \left. U_i^{n} - U_{i-1}^n \right. \right.+ \\[1ex] 
\nonumber
\left. 0.5 (1 -  \mathcal{C}_i) ( U_{i+1}^n - U_i^n) -  0.5 (1 -  \mathcal{C}_{i-1}) ( U_{i}^n - U_{i-1}^n)\right)
\end{eqnarray}
and for $\kappa = 0$ it turns to
\begin{eqnarray}
\label{fexplF}
U_i^{n+1} = U_i^{n} -  \mathcal{C}_i \left( \left. U_i^{n} - U_{i-1}^n \right. \right.+ \\[1ex] 
\nonumber
\left. 0.25 (1 -  \mathcal{C}_i) ( U_{i+1}^n - U_{i-1}^n) -  0.25 (1 -  \mathcal{C}_{i-1}) ( U_{i+1}^n - U_{i-1}^n)\right) \,.
\end{eqnarray}
  

Next we present the leading term of truncation error of the scheme (\ref{fexpl}) for the case of constant velocity, i.e. $\mathcal{C}_i \equiv \mathcal{C}$, that takes the form
\begin{equation}
 {h^3}/{12} \,  (1-|\mathcal{C}|) \left([\mathcal{C}]^+\left( 1 - 2 \mathcal{C} - 3 \kappa  \right) + [\mathcal{C}]^- \left( 1 + 2 \mathcal{C}+ 3 \kappa  \right) \right)
\partial_{xxx} u_i^n \,.
\nonumber
\end{equation}
Therefore, the choice
\begin{equation}
\label{ex3rd}
\kappa = \hbox{sign}(\mathcal{C}) ( 1 -2 |\mathcal{C}| ) / 3
\end{equation}
gives the $3^{rd}$ order accurate scheme in the case of constant velocity.
We search for similar property later in the case of fully implicit and semi-implicit schemes.

To derive numerical stability condition for (\ref{fexpl}) we consider the case of constant velocity.
The numerical stability condition $|S| \le 1$ is obtained for $\kappa \in [-1,1]$ and $\mathcal{C} \in [-1,1]$ that is in the agreement with available theoretical results of von Neumann stability analysis in e.g. \cite{wes01}. It is interesting to note that for $\kappa=-\hbox{sign}(\mathcal{C})$ the numerical stability condition  is obtained for $\mathcal{C} \in [-2,2]$.

We investigate now if a $2^{nd}$ order accurate fully implicit $\kappa$-scheme to solve (\ref{adveq1d}) can be derived analogously to (\ref{fexpl})  with less restrictive numerical stability condition.
It is now rather straightforward to do. Firstly, instead of (\ref{ext}) we use the Taylor series in a backward manner
\begin{equation}
\label{imt}
u_i^{n} = u_i^{n+1} + \sum_{m=1}^p \frac{(-1)^m}{m!} \tau^m \partial^m_t u_i^{n+1}  + \mathcal{O}(\tau^{p+1}) \,.
\end{equation}
Applying the Cauchy-Kowalewski procedure (\ref{ck01}) - (\ref{ck11}) at $t^{n+1}$ with (\ref{imt}) for $p=2$ we obtain
\begin{eqnarray}
\label{accur2i}
u_i^{n} = u_i^{n+1} + \tau V_i \partial_x u_i^{n+1} + 0.5 \tau^2 V_i \partial_{x} (V \partial_x u)_i^{n+1}
+ \mathcal{O}(\tau^3) \,.
\end{eqnarray}
Using appropriate approximations analogous to (\ref{upwgrad1d}) - (\ref{dxx}) in (\ref{accur2i}), the fully implicit $\kappa$-scheme is obtained
\begin{eqnarray}
\label{impl1d}
U_i^{n+1} +  \tau [V_i]^+ \partial_x^- \left( U_i^{n+1} + 0.5 (h + \tau  V_i) \partial_x^{\kappa} U_i^{n+1})\right) + \\[1ex] 
\nonumber
+ \left. \tau [V_i]^- \partial_x^+ \left( U_i^{n+1} - 0.5 (h - \tau  V_i) \partial_x^{\kappa} U_i^{n+1})\right) \right.  = U_i^{n}\,.
\end{eqnarray}
Similarly to (\ref{ex3rd}), one can prove that the choice
\begin{equation}
\label{im3rd}
\kappa = \hbox{sign}(\mathcal{C}) ( 1 + 2 |\mathcal{C}| ) / 3
\end{equation}
gives the $3^{rd}$ order accurate scheme in the case of constant velocity.

We present now the numerical stability conditions for (\ref{impl1d})  for the case of constant velocity. 
The numerical von Neumann stability analysis suggests that for $\kappa \le 0$ the fully implicit $\kappa$-scheme is unconditionally stable for $\mathcal{C}\ge 0$. In the case $\mathcal{C} \le 0$ the unconditional numerical stability is indicated for $\kappa \ge 0$.
These conditions can be seen as an advantage when compared to fully explicit $\kappa$-schemes. The price to pay is that a system of linear algebraic equations has to be solved in each time step to obtain the values $U_i^{n+1}$. 

Unfortunately, the other interesting choices of $\kappa$ give only restrictive numerical stability conditions. For instance 
the value $\kappa=1/3$ gives the conditional numerical stability for $0 \le \mathcal{C} \le 2$, the $3^{rd}$ order accurate scheme (\ref{im3rd}) gives it only for $|\mathcal{C}| \le 0.5$. Moreover, the choice $\kappa=\hbox{sign}(V)$ gives unstable numerical scheme for $\mathcal{C} \in [-1,1]$. 

Now we are ready to present the semi-implicit variant of $\kappa$-scheme.
The main idea is to apply the partial Cauchy-Kowalewski procedure
\begin{eqnarray}
\label{ck01i}
\partial_t u_i^{n+1} = - V_i \partial_x u_i^{n+1} \,,\\
\label{ck11i}
\partial_{tt} u_i^{n+1} = - V_i \partial_{tx} u_i^{n+1} 
\end{eqnarray}
and to avoid the full procedure by skipping the replacement of $\partial_{tx} u$ in (\ref{ck11i}). Using (\ref{ck01i}) and (\ref{ck11i}) with (\ref{imt}) for $p=2$ we obtain
\begin{eqnarray}
\label{accur4}
u_i^{n} =  u_i^{n+1} + \tau V_i \partial_x u_i^{n+1} - 0.5 \tau^2 V_i \partial_{t x}  u_i^{n+1}
+ \mathcal{O}(\tau^3) \,.
\end{eqnarray}
Now using the approximation (\ref{upwgrad1d}) in (\ref{accur4})  for the term containing $\partial_x u_i^{n+1}$ and for the other term the following approximation
$$
\tau \, \partial_{tx} u_i^{n+1} \approx \tau \, \partial_t^- \partial_x^{\kappa} U_i^{n+1} = \partial_x^{\kappa} U_i^{n+1} - \partial_x^{\kappa} U_i^{n} \,,
$$
one obtains 
\begin{eqnarray}
\nonumber
U_i^{n} = U_i^{n+1} +  \tau [V_i]^+ \left(\partial_x^- U_i^{n+1} + 0.5 h \partial_x^-  \partial_x^{\kappa} U_{i}^{n+1} \right) + \\[1ex]
\nonumber
+ \left. \tau [V_i]^- \left(\partial_x^+ U_i^{n+1} - 0.5 h \partial_x^+ \partial_x^{\kappa} U_{i}^{n+1}  \right) \right. 
-0.5 \tau V_i \left( \partial_x^{\kappa} U_i^{n+1} - \partial_x^{\kappa} U_i^{n} \right)\,.
\end{eqnarray}
After simple algebraic manipulations the semi-implicit $\kappa$-scheme can be written in the form of (\ref{all1d})
\begin{eqnarray}
\label{si1d}
U_i^{n+1} +  \tau V_i \left(\partial_x^{\mp} U_i^{n+1} - 0.5 \partial_x^{\kappa} U_{i\mp 1}^{n+1} \right) 
= \left.
U_i^{n} -   0.5 \tau V_i \partial_x^{\kappa} U_{i}^{n} \right. \,,
\end{eqnarray}
where one has to replace $\mp$ in $\partial_x^{\mp}$ and in $i\mp 1$ with opposite signs with respect to $V_i$, i.e. with $-$ if $V_{i}>0$ and with $+$ if $V_{i}<0$.

Looking at the truncation error of  (\ref{si1d}) the scheme is $2^{nd}$ order accurate for variable velocity case and for arbitrary value of $\kappa$. The variable choice $\kappa$ with respect to $i$
\begin{equation}
\label{kappasi}
\kappa = \hbox{sign}(\mathcal{C}_i) (1- |\mathcal{C}_i|)/3
\end{equation}
gives in the case of constant velocity the $3^{rd}$ order accurate scheme. As we show later when considering the scheme in several dimensions, the choice (\ref{kappasi}) can be advantageous also for variable velocity case.

To use (\ref{si1d}) one has to solve a linear system of algebraic equations. If $V(x)$ does not change its sign, the matrix of resulting system has a three-diagonal form with two off-diagonals strictly either below or up to the main diagonal for all values of $\kappa$. Consequently, the system can be solved in one step using a forward or backward substitution.

The most important property of semi-implicit $\kappa$-scheme is concerning its numerical stability condition. The scheme (\ref{si1d}) exploits only the single value $V_i$, so the stability analysis is valid without assumption of constant velocity. The numerical von Neumann stability analysis indicates unconditional stability for arbitrary $\kappa \le 1$ if $V_i>0$ and for $\kappa \ge -1$ if $V_i<0$.

Analogously to (\ref{fexplLW}) and (\ref{fexplF}) we present two particular variants of (\ref{si1d}) for the case of variable velocity to present them in a more clear way.
Particularly for $\kappa = \hbox{sign}(V_i)$ the scheme (\ref{si1d}) takes the form
\begin{eqnarray}
\label{siLW}
U_i^{n+1} +  0.5 \, |\mathcal{C}_i| ( U_i^{n+1} - U_{i\mp 1}^{n+1} ) = U_i^{n} -  0.5 \, |\mathcal{C}_i| ( U_{i\pm 1}^{n} - U_{i}^n ) \,,
\end{eqnarray}
where the signs in $\mp$  are chosen opposite and the signs in $\pm$ identical to the sign of $V_i$. The scheme (\ref{siLW}) is introduced in \cite{mo11} in a finite volume context, see also \cite{mo10,mou14}. Furthermore for $\kappa \equiv 0$ the scheme (\ref{si1d}) turns to
\begin{eqnarray}
\label{siF}
U_i^{n+1} +  0.25 \, |\mathcal{C}_i| ( 3 U_i^{n+1} - 4 U_{i\mp 1}^{n+1} +U_{i\mp 2}^{n+1} ) = U_i^{n} -  0.25 \, |\mathcal{C}_i| ( U_{i+1}^{n} - U_{i-1}^n ) \,.
\end{eqnarray}
The scheme (\ref{siF}) is introduced in \cite{fmu14} in the finite volume context.

\begin{remark}
It is clear that the partial Cauchy-Kowalewski procedure can be applied not only to (\ref{imt}), but also to the Taylor series (\ref{ext}). Analogously to (\ref{si1d}) one obtain in this case
\begin{eqnarray}
\label{exception1db}
U_i^{n+1} +  0.5 \, \tau V_i \partial_x^{\kappa} U_{i}^{n+1} 
= U_i^{n} - \tau V_i (\partial_x^{\mp} U_i^{n} - 0.5  \partial_x^{\kappa} U_{i\mp 1}^{n} ) 
\end{eqnarray}
where again the sign in $\mp$ shall be chosen opposite to the sign of $V_i$.

The numerical stability condition of (\ref{exception1db}) is obtained for $\kappa \le 1$ in the form $\mathcal{C}_i \in [0,1]$ and  for $\kappa \ge -1$ in the form for $\mathcal{C}_i \in [-1,0]$. For the case $\kappa = \hbox{sign}(V_i)$ one has that $|S|=1$ so such variant of (\ref{exception1db}) has unconditional numerical stability.

\end{remark}


\section{Two-dimensional case}
\label{sec-2d}

The representative 2D advection equation takes the form
\begin{equation}
\partial_t u(x,y,t) + \vec{V}(x,y) \cdot \nabla u(x,y,t) = 0 \,, \quad
u(x,y,0)=u^0(x,y) \,,
\label{adveq}
\end{equation}
where $\vec{V}=(V(x,y),W(x,y))$.
We restrict ourselves to Cartesian grids that can be obtained from uniform one-dimensional grids using standard dimension by dimension extension. We denote the uniform space discretization step by $h$.  We aim to suggest numerical schemes to determine the approximate values $U_{ij}^n \approx u_{i j}^{n}$ where  $u_{ij}^n := u(x_i,y_j,t^n)$.

The extension of general 1D scheme (\ref{all1d}) for two-dimensional case can be written in the form
\begin{eqnarray}
\label{all2d}
U_{ij}^{n+1} + \\
\nonumber
\sum_{k=-2}^{2} \left( \alpha^x_{i j k} U_{i+k j}^{n+1} + \alpha^y_{i j k} U_{i j+k}^{n+1} \right) = 
U_{ij}^{n} +\sum_{k=-2}^{2} \left( \beta^x_{i j k} U_{i+k j}^{n} + \beta^y_{i j k} U_{i j+k}^{n} \right) ,
\end{eqnarray}
where we adopt the notation of superscripts $x$ and $y$ to relate the coefficients to particular space variable.

Analogously to (\ref{kappagrad}), the following approximation of gradients $( \partial_x u_{i j}^n,\partial_y u_{i j}^n)$ is used,
\begin{eqnarray}
\label{grad2d}
2 h \, \partial_x^{\kappa} U_{i j}^n =  (1-\kappa^x) \partial_x^- U_{i j}^n + (1+\kappa^x) \partial_x^+ U_{i j}^n \,, \\
\label{grad2dB}
2 h \, \partial_y^{\kappa} U_{i j}^n =  (1-\kappa^y) \partial_y^- U_{i j}^n + (1+\kappa^y) \partial_y^+ U_{i j}^n\,,  
\end{eqnarray}
where $\partial_x^{-}$, $\partial_x^{+}$, $\partial_y^{-}$, and $\partial_y^{+}$ denote the standard finite differences analogously to (\ref{kappagrad}), and the parameters $\kappa^x$ and $\kappa^y$ are free to choose. 

To provide numerical stability analysis of (\ref{all2d}), one extends the one-dimensional treatment (\ref{eps1d}) - (\ref{allstab}) by using a grid function $\epsilon_{ij}^n=\epsilon(x_i,y_j,t^n)$ 
when the amplification factor $S$ takes the form analogously to (\ref{allstab})
\begin{eqnarray}
\label{allstab2d}
S = \frac{
1 + \sum_{k=-2}^{2} (\beta^x_{i j k} \exp(\imath k x)+\beta^y_{i j k} \exp(\imath k y)) 
}{
1+ \sum_{k=-2}^{2} \alpha^x_{i j k}  \exp(\imath k x) + \sum_{k=-2}^{2} \alpha^y_{i j k}  \exp(\imath k y) 
}
\end{eqnarray}
and $x,y \in (-\pi,\pi)$.

Before deriving particular numerical schemes of the form (\ref{all2d}) we note that the Cauchy-Kowalewski procedure takes now more involved form than (\ref{ck01}) - (\ref{ck11}) in 1D case, namely
\begin{eqnarray}
\label{taylor2dA}
\partial_{t} u_{ij}^n = -V_{ij} \partial_{x} u_{ij}^n - W_{ij} \partial_{y} u_{ij}^n \,,  \,\, \partial_{tt} u_{ij}^n = -V_{ij} \partial_{tx} u_{ij}^n - W_{ij} \partial_{ty} u_{ij}^n \,,  \,\, \\[1ex]
\label{taylor2dB}
\partial_{tx} u_{ij}^n = -\partial_x (V \partial_{x} u)_{ij}^n - \partial_x (W \partial_{y} u)_{ij}^n,  
\partial_{ty} u_{ij}^n = -\partial_y (V \partial_{x} u)_{ij}^n - \partial_y (W \partial_{y} u)_{ij}^n.
\end{eqnarray}

Clearly, when using the full Cauchy-Kowalewski procedure to derive a fully explicit or fully implicit variant of (\ref{all2d}), one has to approximate, due to (\ref{taylor2dB}), the mixed spatial derivative of $u$. This can not be done using the stencil prescribed by (\ref{all2d}), consequently the dimension by dimension extension of 1D schemes (\ref{fexpl}) and (\ref{impl1d}) in the form of (\ref{all2d}) are not $2^{nd}$ order accurate schemes as already well-known from literature for fully explicit variants, see e.g. \cite{bil97b,lev02}.

On the other hand, when using the partial Cauchy-Kowalewski procedure (\ref{taylor2dA}) without (\ref{taylor2dB}) no mixed spatial derivatives are involved. Consequently, it can be shown that the dimension by dimension extension of semi-implicit $\kappa$-scheme (\ref{si1d}) is $2^{nd}$ order accurate for variable velocity $\vec{V}(x,y)$ and for arbitrary values of $\kappa^x$ and $\kappa^y$ in (\ref{grad2d}). 

Similarly to (\ref{si1d}) we write the semi-implicit $\kappa$-scheme in the compact (upwind) way 
\begin{eqnarray}
\nonumber
U_{i j}^{n+1} + 
\tau V_{i j} \left(\partial_x^{\mp} U_{i j}^{n+1} - 0.5 \partial_x^{\kappa} U_{i\mp 1\, j}^{n+1} \right) + 
\tau W_{i j} \left(\partial_y^{\mp} U_{i j}^{n+1} - 0.5 \partial_y^{\kappa} U_{i\, j\mp 1}^{n+1} \right) 
= \\[1ex]
\label{si2d}
U_i^{n} - 0.5 \tau \left(V_{i j} \partial_x^{\kappa} U_{i j}^{n} + W_{i j} \partial_y^{\kappa} U_{i j}^{n} \right) ,
\end{eqnarray}
where one has to replace $\mp$ in $\partial_x^{\mp}$ and in $i\mp 1$ with $-$ if $V_{ij}>0$ and with $+$ if $V_{ij}<0$, compare also with (\ref{si1d}),  and analogously for the cases related to the sign of $W_{ij}$. Note that the scheme (\ref{si2d}) can be extended to higher dimensional case.

It appears that, opposite to one-dimensional case, the semi-implicit $\kappa$-scheme (\ref{si2d}) is not numerically unconditional stable in general. It is not easy to characterize the numerical stability conditions for all cases as the amplification factor $S$ in (\ref{allstab2d}) depends on six free parameters: $x$, $y$, $\kappa^x$, $\kappa^y$ and two (directional and signed) grid Courant numbers
\begin{eqnarray}
\label{cn2d}
\mathcal{C}_{ij}= \frac{\tau V_{ij}}{h} \,, \,\, \mathcal{D}_{i j} = \frac{\tau W_{ij}}{h} \,.
\end{eqnarray}

We found that in general the numerical stability condition for (\ref{si2d}) is indicated for arbitrary $|\kappa^x| \le 1$ and $|\kappa^y| \le 1$  if $|\mathcal{C}_{ij}| \le 4$ and $|\mathcal{D}_{ij}| \le 4$. In section on numerical experiments we report unstable numerical solutions for parameters outside of these intervals.

Furthermore we give more details on stability for the variants of (\ref{si2d}) that have been published elsewhere and that are used in section on numerical experiments.

The choice $\kappa^x=\hbox{sign}(V_{ij}), \kappa^y=\hbox{sign}(W_{ij})$ in (\ref{si2d}) gives the so-called IIOE scheme (Inflow Implicit / Outflow Explicit) published in a finite volume form in \cite{mou14}. The scheme gives $|S|=1$ for all values of $\mathcal{C}_{ij},\mathcal{D}_{ij},x,y$. Consequently it can be considered as unconditionally stable, but it does not damp possible oscillations in numerical solution, see also corresponding results in numerical experiments. 

The choice $\kappa^x= \kappa^y=0$ in (\ref{si2d}) is used in a finite volume form in \cite{fmu14} and it gives numerical stability condition $|S| \le 1$ for $|\mathcal{C}_{ij}| \le 7.396$ and $|\mathcal{D}_{ij}| \le 7.396$, but the value $|S|$ can be larger than 1 otherwise. For instance the maximal value of $|S|$ is around $1.00013$ and $1.04538$ for maximal Courant numbers $8$ and $16$, respectively. The code in Mathematica \cite{mat16} by which these numerical stability results are obtained is available by a request.

Fortunately one can extend the semi-implicit $\kappa$-scheme in two-dimensional case to such a form for which unconditional numerical stability is indicated.

\subsection{Corner Transport Upwind extension}

In what follows we apply the idea of so-called Corner Transport Upwind (CTU) scheme \cite{col90}. We follow \cite{lev02} where it is used to extend the fully explicit schemes of the form (\ref{all2d}) to take into account the mixed derivatives in (\ref{taylor2dB}). 

Firstly, we analyze the truncation error of (\ref{si2d}) for the choice of $\kappa$ parameters analogous to (\ref{kappasi}) for the case of constant velocity $\vec{V}$. As this choice will be used in numerical schemes also for variable velocity case, we keep the indexing in the notation $(V_{ij},W_{ij})$ of discrete velocity values.

We consider the Taylor series
\begin{equation}
\label{imt2d}
u_{ij}^{n} = u_{ij}^{n+1} + \sum_{m=1}^p \frac{(-1)^m}{m!} \tau^m \partial^m_t u_{ij}^{n+1}  + \mathcal{O}(\tau^{p+1}) \,,
\end{equation}
for $p=3$ and replace the first and second time derivatives in (\ref{imt2d}) using (\ref{taylor2dA}) at $t^{n+1}$. The third term can be replaced for constant velocity $(V_{ij},W_{ij})$ by
\begin{eqnarray}
\label{taylor2d3A}
\partial_{ttt} u_{ij}^{n+1} = -V_{ij} \partial_{ttx} u_{ij}^{n+1} - W_{ij} \partial_{tty} u_{ij}^{n+1} \,.
\end{eqnarray}
In our approach we have to furthermore substitute to (\ref{taylor2d3A}) the relations valid for the constant vector $(V_{ij},W_{ij})$
\begin{eqnarray}
\label{taylor2d3B}
\partial_{ttx} u_{ij}^{n+1} = -V_{ij} \partial_{txx} u_{ij}^{n+1} - W_{ij}  \partial_{txy} u_{ij}^{n+1} \,,  \,\,\\[1ex]
\label{taylor2d3C}
\partial_{tty} u_{ij}^{n+1} = -V_{ij} \partial_{txy} u_{ij}^{n+1} - W_{ij} \partial_{tyy} u_{ij}^{n+1}\,.
\end{eqnarray}

One can show that for the choices analogous to (\ref{kappasi})
\begin{equation}
\kappa^x = \hbox{sign}(\mathcal{C}_{i j}) (1- |\mathcal{C}_{i j}|)/3 \,, \,\,
\kappa^y = \hbox{sign}(\mathcal{D}_{i j}) (1- |\mathcal{D}_{i j}|)/3 \,, 
\label{kappasi2d}
\end{equation}
the spatial derivatives $\partial_{xxx}u_{i j}^{n+1}$ and $\partial_{yyy} u_{i j}^{n+1}$  are canceled in the truncation error of (\ref{si2d}) analogously to 1D case. Nevertheless, the following third order derivatives term 
will remain
\begin{equation}
\frac{\tau^2}{12} V_{i j} W_{i j} \left(2 \tau \partial_{txy} u_{i j}^{n+1} - h \partial_{xxx} u_{i j}^{n+1} - h \partial_{yyy} u_{i j}^{n+1} \right) \,.
\label{mixed}
\end{equation}
Consequently, the scheme (\ref{si2d}) with (\ref{kappasi2d}) can not be $3^{rd}$ order accurate in the case of constant velocity.

A relatively simple form of (\ref{mixed}) motivates us to extend the scheme (\ref{si2d}) using the approach of Corner Transport Upwind (CTU) scheme \cite{col90,lev02}. This extension consists of adding additional discretization terms to (\ref{si2d}) that contain, additionally to (\ref{all2d}), also the corner (diagonal) values of numerical solution. It can be done in such a way that the $2^{nd}$ order accuracy is preserved and the scheme becomes $3^{rd}$ order accurate in the case of constant velocity $\vec{V}$. 

In fact one can derive two such schemes that differ only in the explicit part.
Firstly,
\begin{eqnarray}
\nonumber
U_{i j}^{n+1} + 
\tau V_{i j} \left(\partial_x^{\mp} U_{i j}^{n+1} - 0.5 \partial_x^{\kappa} U_{i\mp 1\, j}^{n+1} \right) + 
\tau W_{i j} \left(\partial_y^{\mp} U_{i j}^{n+1} - 0.5 \partial_y^{\kappa} U_{i\, j\mp 1}^{n+1} \right) + \\[1ex]
\nonumber
|\mathcal{C}_{ij} \mathcal{D}_{ij}|/6 \left( U_{i j}^{n+1}+U_{i\mp1\, j\mp1}^{n+1}-U_{i\mp1\, j}^{n+1}-U_{i\, j\mp1}^{n+1} \right)
= \\[1ex]
\nonumber
U_i^{n} - 0.5 \tau \left(V_{i j} \partial_x^{\kappa} U_{i j}^{n} + W_{i j} \partial_y^{\kappa} U_{i j}^{n} \right) + \\[1ex]
\label{ctusi2d1}
|\mathcal{C}_{ij} \mathcal{D}_{ij}|/12 \left( 2 U_{i j}^n + U_{i\pm1 j\pm1}^{n}+U_{i\mp1\, j\mp1}^{n}-U_{i+1\, j}^{n}-U_{i\, j+1}^n-U_{i-1\, j}^{n}-U_{i\, j-1}^{n} \right).
\end{eqnarray}
and, secondly,
\begin{eqnarray}
\nonumber
U_{i j}^{n+1} + 
\tau V_{i j} \left(\partial_x^{\mp} U_{i j}^{n+1} - 0.5 \partial_x^{\kappa} U_{i\mp 1\, j}^{n+1} \right) + 
\tau W_{i j} \left(\partial_y^{\mp} U_{i j}^{n+1} - 0.5 \partial_y^{\kappa} U_{i\, j\mp 1}^{n+1} \right) + \\[1ex]
\nonumber
|\mathcal{C}_{ij} \mathcal{D}_{ij}|/6 \left( U_{i j}^{n+1}+U_{i\mp1\, j\mp1}^{n+1}-U_{i\mp1\, j}^{n+1}-U_{i\, j\mp1}^{n+1} \right)
= \\[1ex]
\nonumber
U_i^{n} - 0.5 \tau \left(V_{i j} \partial_x^{\kappa} U_{i j}^{n} + W_{i j} \partial_y^{\kappa} U_{i j}^{n} \right) - \\[1ex]
\label{ctusi2d2}
|\mathcal{C}_{ij} \mathcal{D}_{ij}|/12 \left( 2 U_{i j}^n + U_{i\mp1 j\pm1}^{n}+U_{i\pm1\, j\mp1}^{n}-U_{i+1\, j}^{n}-U_{i\, j+1}^n-U_{i-1\, j}^{n}-U_{i\, j-1}^{n} \right).
\end{eqnarray}
In (\ref{ctusi2d1}) and (\ref{ctusi2d2}) the identical convention is used for $\mp$ and $\pm$ as for (\ref{si2d}). 

Concerning the stability property, the numerical von Neumann stability analysis indicates that the schemes (\ref{ctusi2d1}) and (\ref{ctusi2d2}) are unconditionally stable for variable velocity $\vec{V}_{ij}$ if the $\kappa$ parameters are chosen as in (\ref{kappasi2d}). In fact any convex linear combination of these two schemes appears to have such unconditional numerical stability, the corresponding code in Mathematica \cite{mat16} by which these numerical results are indicated is available by a request.
Moreover, the choice (\ref{kappasi2d}) in (\ref{ctusi2d1}) and (\ref{ctusi2d2}) results in the third order accurate scheme if the velocity vector $\vec{V}$ is constant.

\section{Numerical experiments}
\label{sec-num}
In what follows we illustrate the properties of semi-implicit $\kappa$-schemes for some benchmark examples. Note that in all examples the resulting linear algebraic systems are solved by Gauss-Seidel iterations using a strategy of so-called fast sweeping method \cite{zao04}) where at most 3 sweeps (i.e. 12 Gauss-Seidel iterations) are used.

To illustrate the formal order of accuracy of all methods we start with numerical examples of advection equation with constant velocity $\vec{V}=(0.8,0.9)$. The exact solution is given simply by $u(x,y,t)=u^0(x-0.8 t, y-0.9 t)$, where $u^0=u^0(x,y)$ is a given initial function. To check the implementation of methods we begin with the choice $u^0(x,y)$ being randomly chosen quadratic function.
Using the exact solution to define Dirichlet boundary conditions, we obtain with (\ref{si2d}) for all interesting choices of ${\kappa}^x$ and ${\kappa}^y$ in (\ref{grad2d}) the exact solution up to a machine accuracy for any chosen $N$ and $M$. Choosing as the initial function some cubic polynomial, only the  schemes (\ref{ctusi2d1}) and (\ref{ctusi2d2}) 
give numerical solutions differing from the exact one purely by rounding errors.

Next we compute examples of rotation for several initial profiles $u^0(x,y)$ where the advection velocity is always defined by 
$$
\vec{V}=(-2 \pi y \, , \, 2 \pi x ) \,.
$$
For any initial function $u^0(x,y)$ one has the exact solution given by $u(x,y,t)=u^0(x \cos (2 \pi t) + y \sin(2 \pi t) , y \cos(2 \pi t) - x \sin (2 \pi t))$. In all examples we choose Dirichlet boundary conditions defined by the exact solution. The domain is a square $(-1,1)^2$ and $t \in [0,1]$, so the profile given by $u^0$ rotates once to return to its initial position at $t=1$. 

To estimate the experimental order of convergence (EOC) for all examples, we compute the following error
\begin{equation}
\label{error}
E = h^2 \max_{n=1,..,N} \sum_{i,j=1}^M  |U_{i j }^n - u(x_i,y_j,t^n)| \,.
\end{equation}
We present also the minimal values of numerical solutions to compare it with the exact value that is always equal to $0$. The results are compared for following choices of $\kappa$,
\begin{eqnarray}
\label{kp}
\kappa^x = \hbox{sgn}(\mathcal{C}_{ij}) \,, \quad \kappa^y = \hbox{sgn}(\mathcal{D}_{ij}) \,, \\
\label{km}
\kappa^x = -\hbox{sgn}(\mathcal{C}_{ij}) \,, \quad \kappa^y = - \hbox{sgn}(\mathcal{D}_{ij}) \,, \\
\label{k0}
\kappa^x = \kappa^y = 0 \,,\\
\label{k3}
\kappa^x = \hbox{sgn}(\mathcal{C}_{ij})(1-|\mathcal{C}_{ij})/3 \,, \quad \kappa^y = \hbox{sgn}(\mathcal{D}_{ij})(1-|\mathcal{D}_{ij})/3 \,. 
\end{eqnarray}

%
%

\begin{figure}[h!]
\begin{center}
\includegraphics[width=0.24\columnwidth]{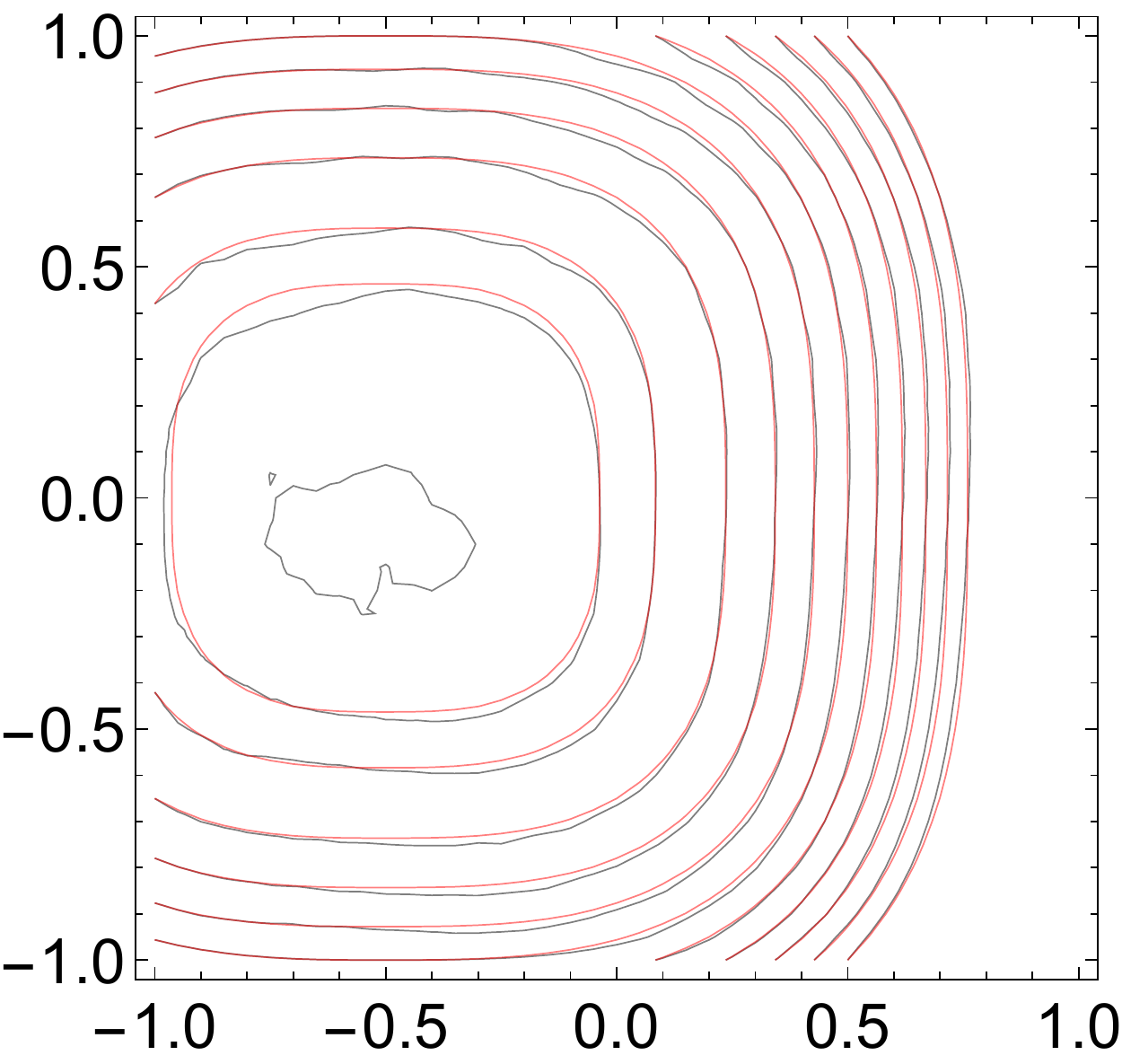}
\includegraphics[width=0.24\columnwidth]{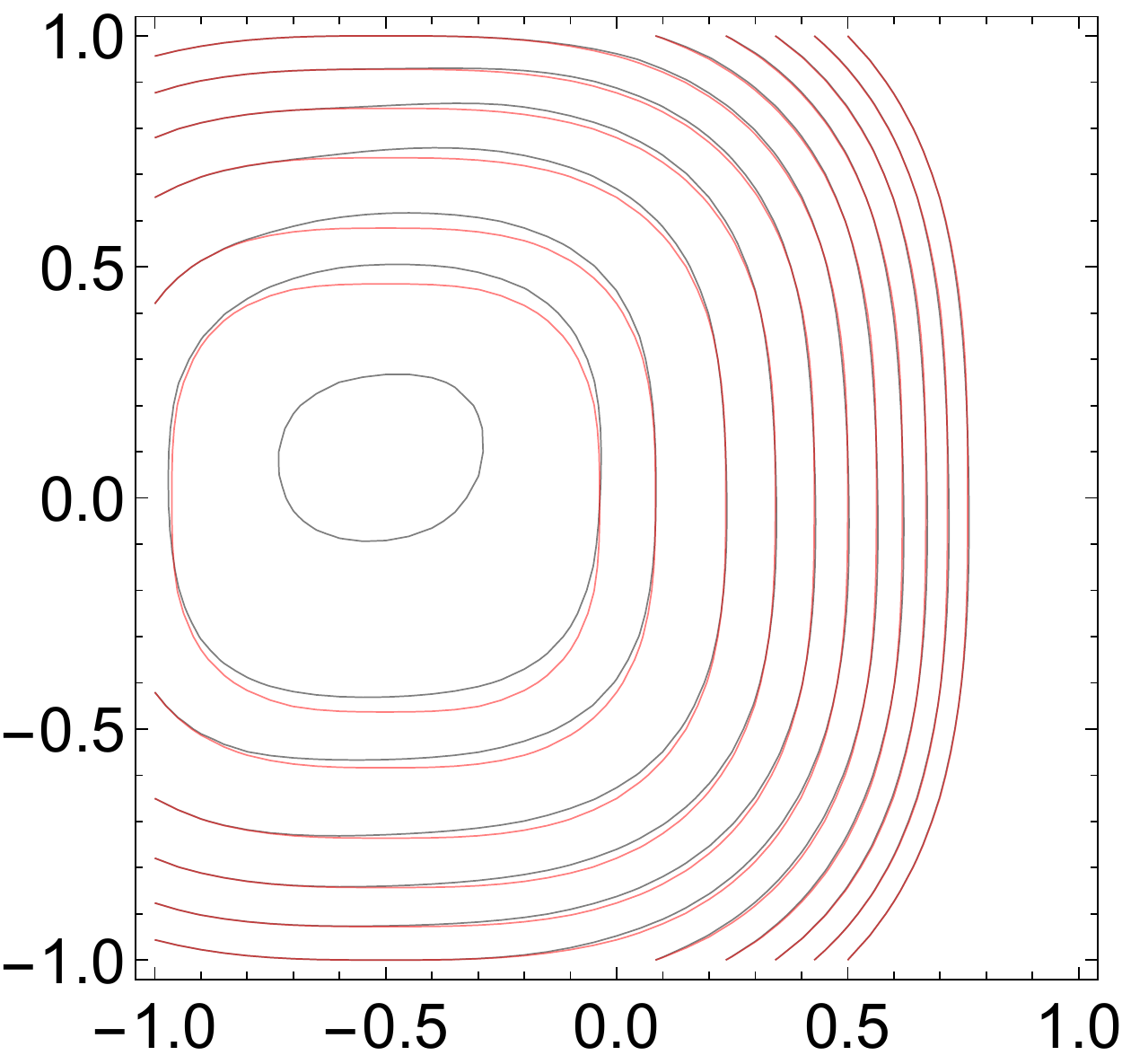}
\includegraphics[width=0.24\columnwidth]{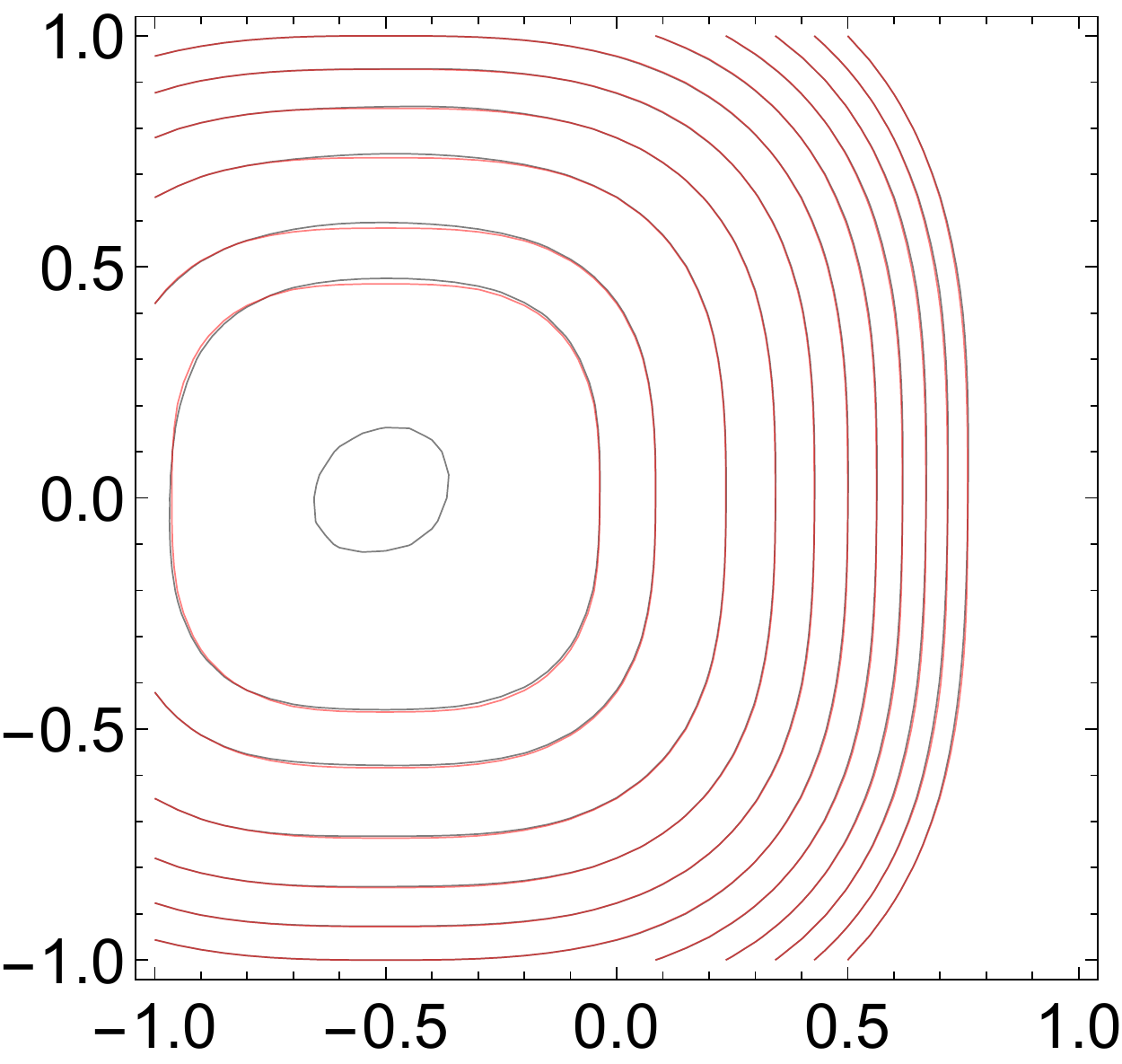}
\includegraphics[width=0.24\columnwidth]{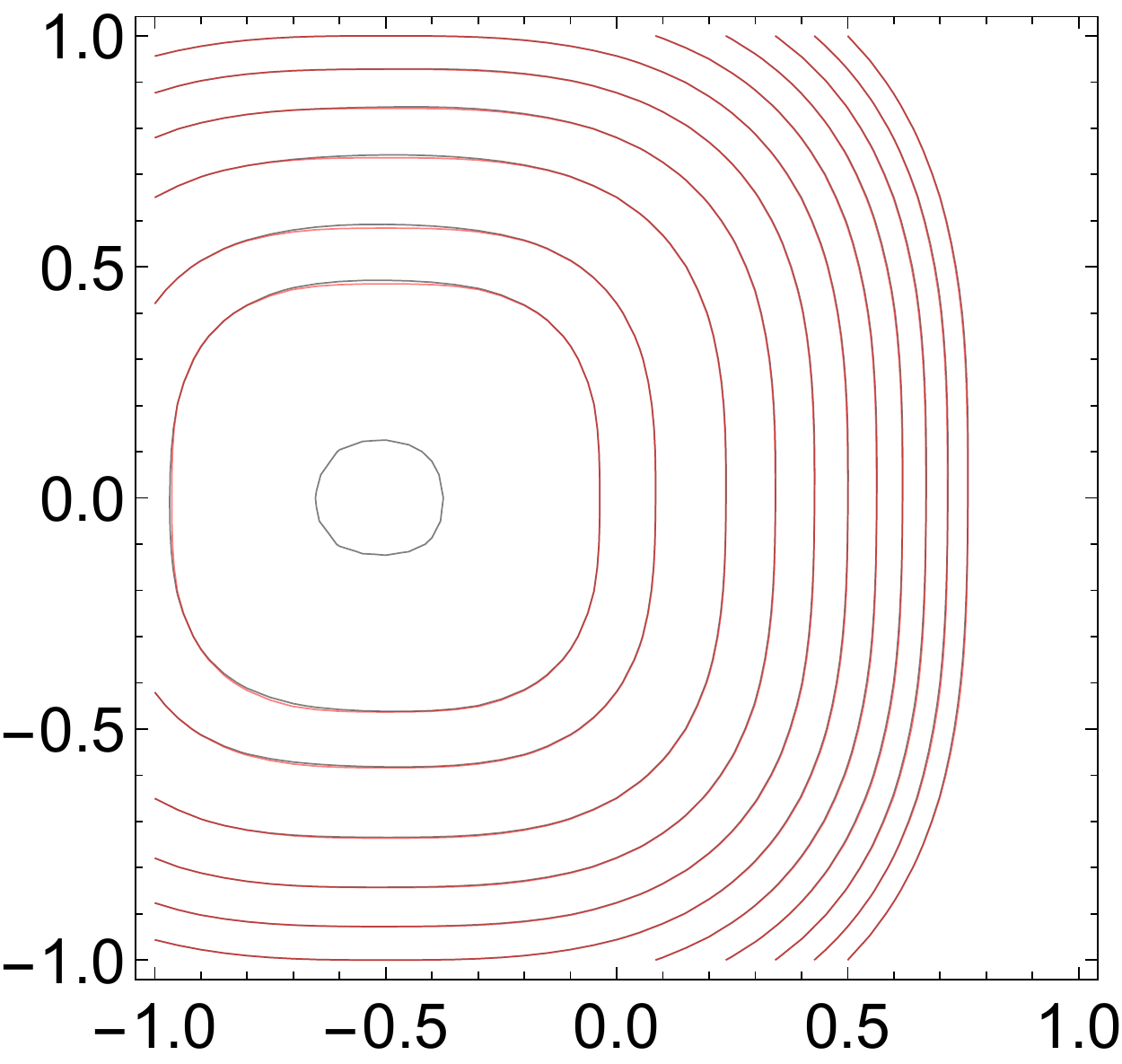}
\caption{Numerical solutions for the rotation of cubic function for the coarse grid $M=40$ and small time step with $N=100$. The red contour lines represent the exact solution, the black contour lines represent the numerical solutions for values $0.0$, $0.2$ up to $2.2$ at $t=1$. The pictures represent the scheme (\ref{si2d}) with (\ref{kp}) - (\ref{k0}) and (\ref{ctusi2d1}) with (\ref{k3}) in that order.
}
\label{fig00}
\end{center}
\end{figure}

Firstly, we choose the cubic function $u^0(x,y)=|x+0.5|(x+0.5)^2 +|y| y^2$. The number of time steps $N=5 M / 2$  is chosen such that the maximum of Courant numbers $|\mathcal{C}_{ij}|$ and $|\mathcal{D}_{ij}|$ is approximately $1.257$.
In Table \ref{tab00} we compare (\ref{si2d}) with all variants (\ref{kp}) - (\ref{k3}) and the CTU extension (\ref{ctusi2d1}) with (\ref{k3}). This order of schemes corresponds to the order from the largest error to the smallest one for this example. Nevertheless for each choice one can observe  the EOC approaching the $2^{nd}$ order accuracy. 
One can observe also small undershootings in each results that diminish quickly with grid refinement. In Figure \ref{fig00} we compare the numerical solutions at $t=1$ with the exact solution. We note that the schemes (\ref{ctusi2d1}) and (\ref{ctusi2d2}) with (\ref{k3}) produces similar results for all presented examples, therefore we present the results  only for the first one. 

\begin{table}[h!]
\begin{center}
\begin{tabular}{||c||c|c||c|c||c|c||c|c||c|c||}
\hline
$M$ & $E$ & $\min$  & $E$ & $\min$ & $E$ & $\min$ & $E$ & $\min$ & $E$ & $\min$ \\
\hline
40  & 54.1 & -15. & 41.2 & -15.  & 13.7 & -4.1 & 11.8 & -3.6 & 9.08 & -3.5 \\
80  & 15.1 & -3.2 & 11.6 & -3.5 & 3.63 & -.86 & 3.12 & -.74 & 2.29 & -.73 \\
160& 3.98 & -.78 & 3.08 & -.82 & .930 & -.19 & .791 & -.16 & .568 & -.15 \\
\hline
\end{tabular}
\caption{The error (\ref{error}) (multiplied by $10^3$) and the minimum of numerical solutions (multiplied by $10^3$) for the rotation of cubic function using (\ref{si2d}) with (\ref{kp})  (the $2^{nd}$ - $3^{rd}$ columns), (\ref{km}) (the $4^{th}$ - $5^{th}$ ones), (\ref{k0}) (the $6^{th}$ - $7^{th}$ ones), (\ref{k3}) (the $8^{th}$ - $9^{th}$ ones) and (\ref{ctusi2d1}) with  (\ref{k3}).  Note that $N=5 M/2$.}
\label{tab00}
\end{center}
\end{table}

\begin{figure}[h!]
\begin{center}
\includegraphics[width=0.24\columnwidth]{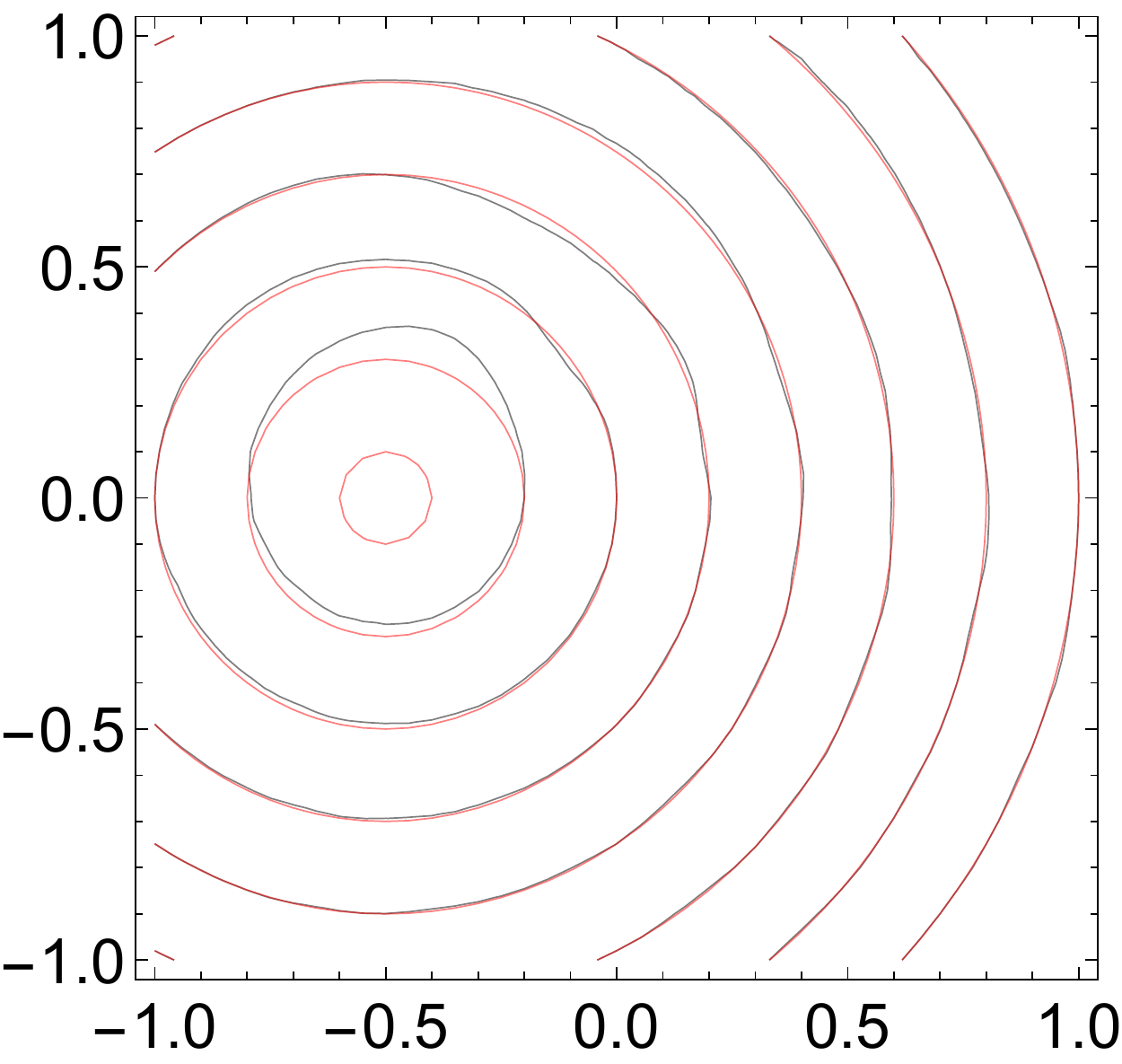}
\includegraphics[width=0.24\columnwidth]{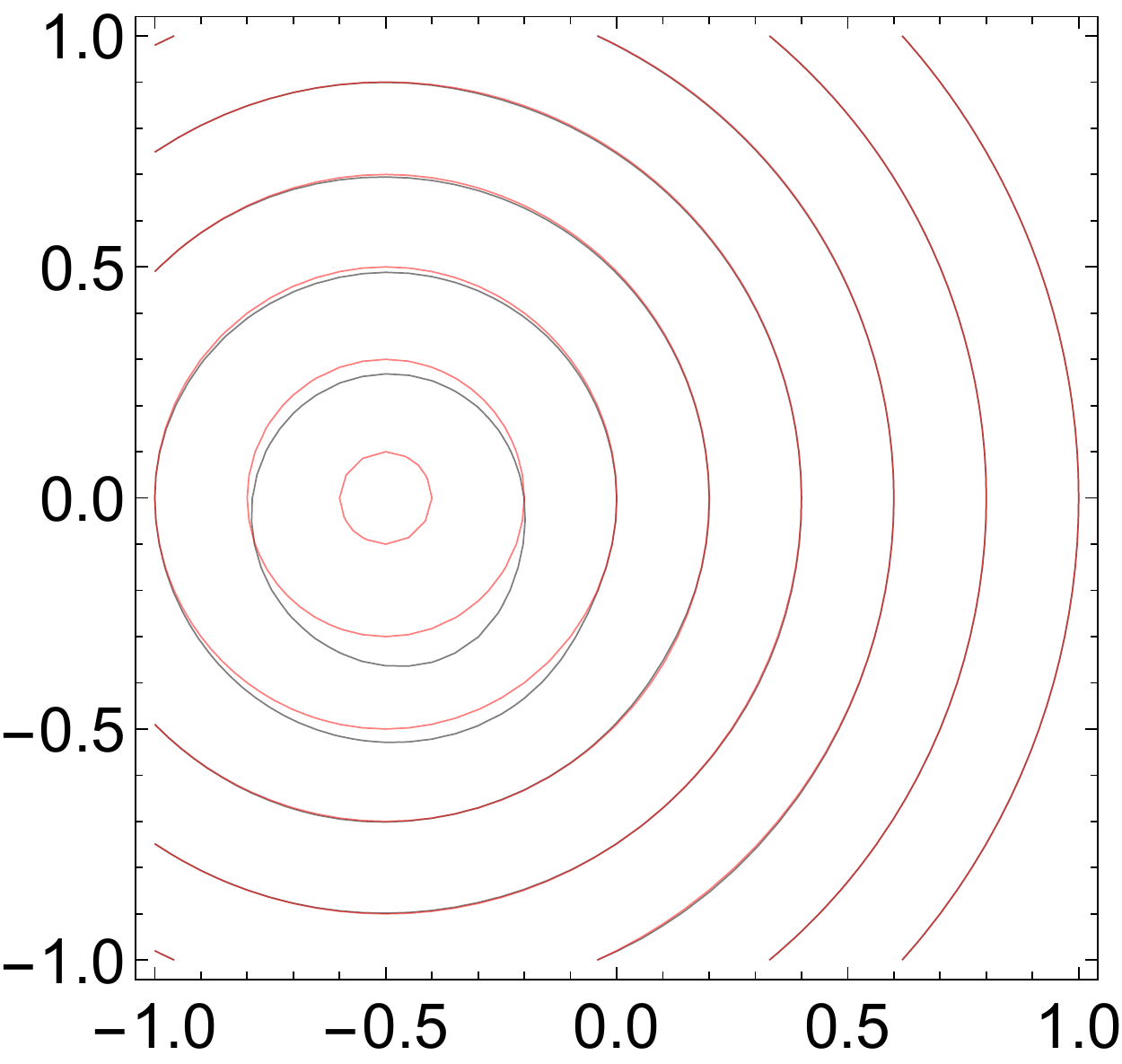}
\includegraphics[width=0.24\columnwidth]{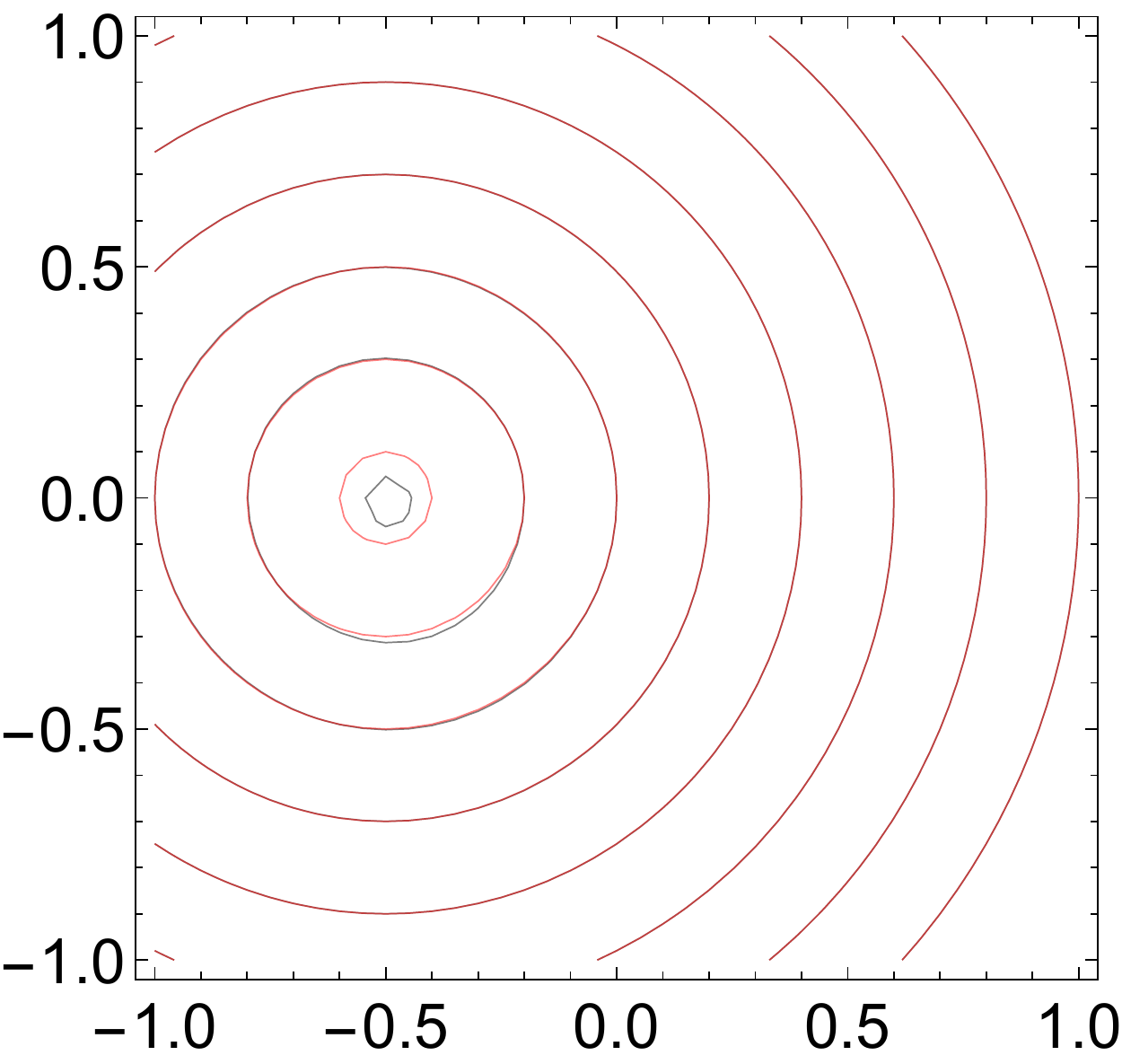}
\includegraphics[width=0.24\columnwidth]{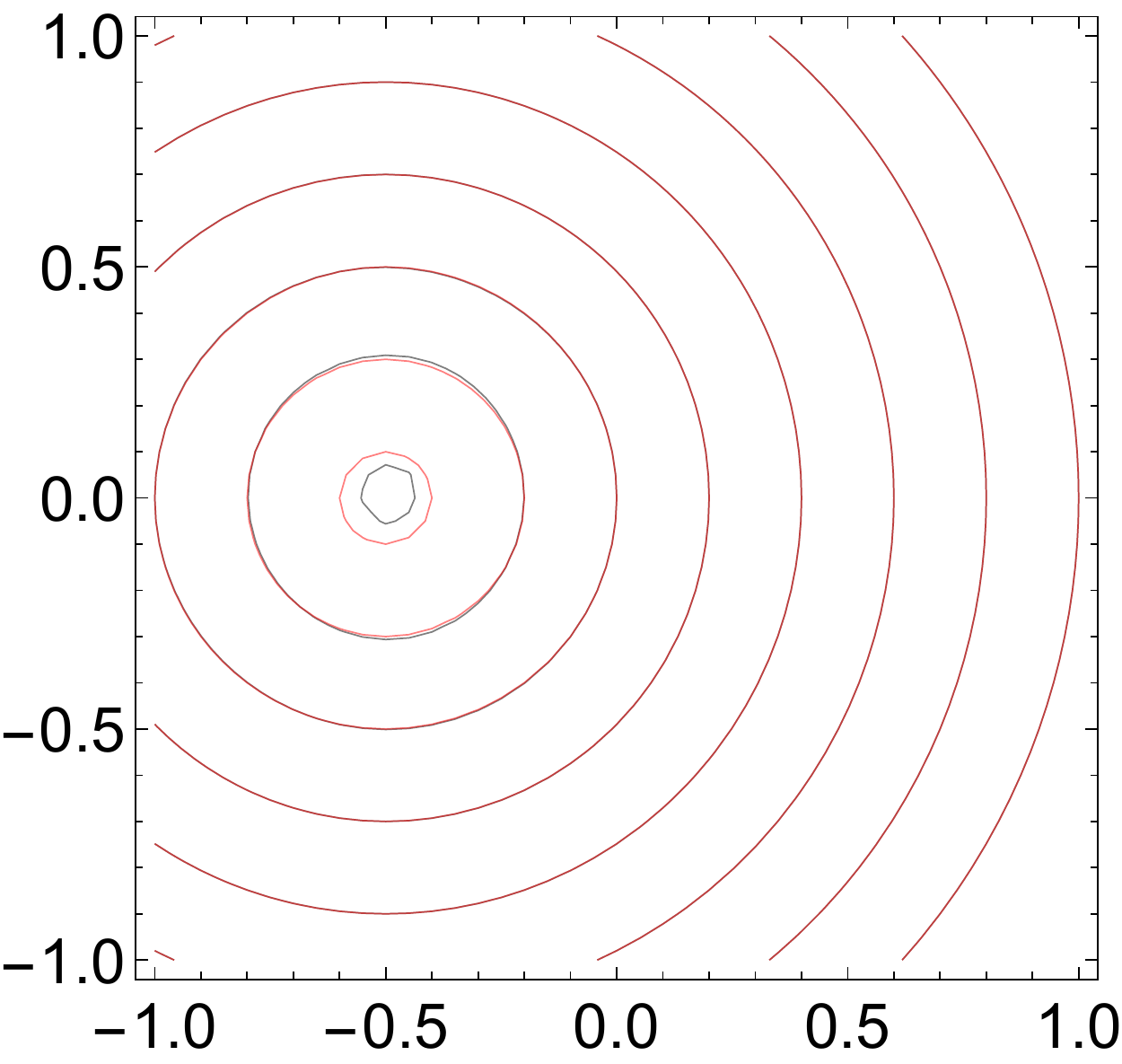}\\[1ex]
\includegraphics[width=0.24\columnwidth]{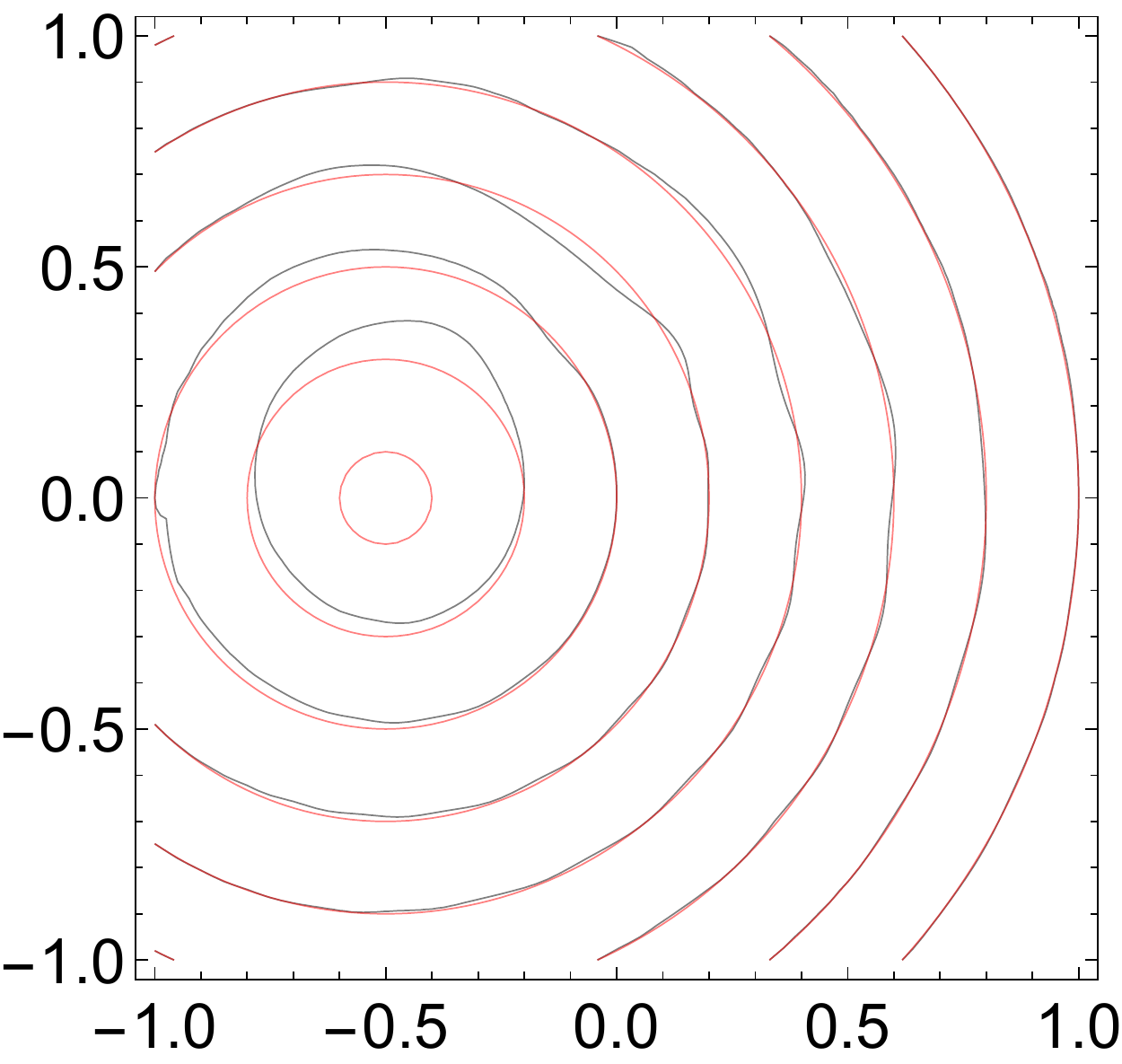}
\includegraphics[width=0.24\columnwidth]{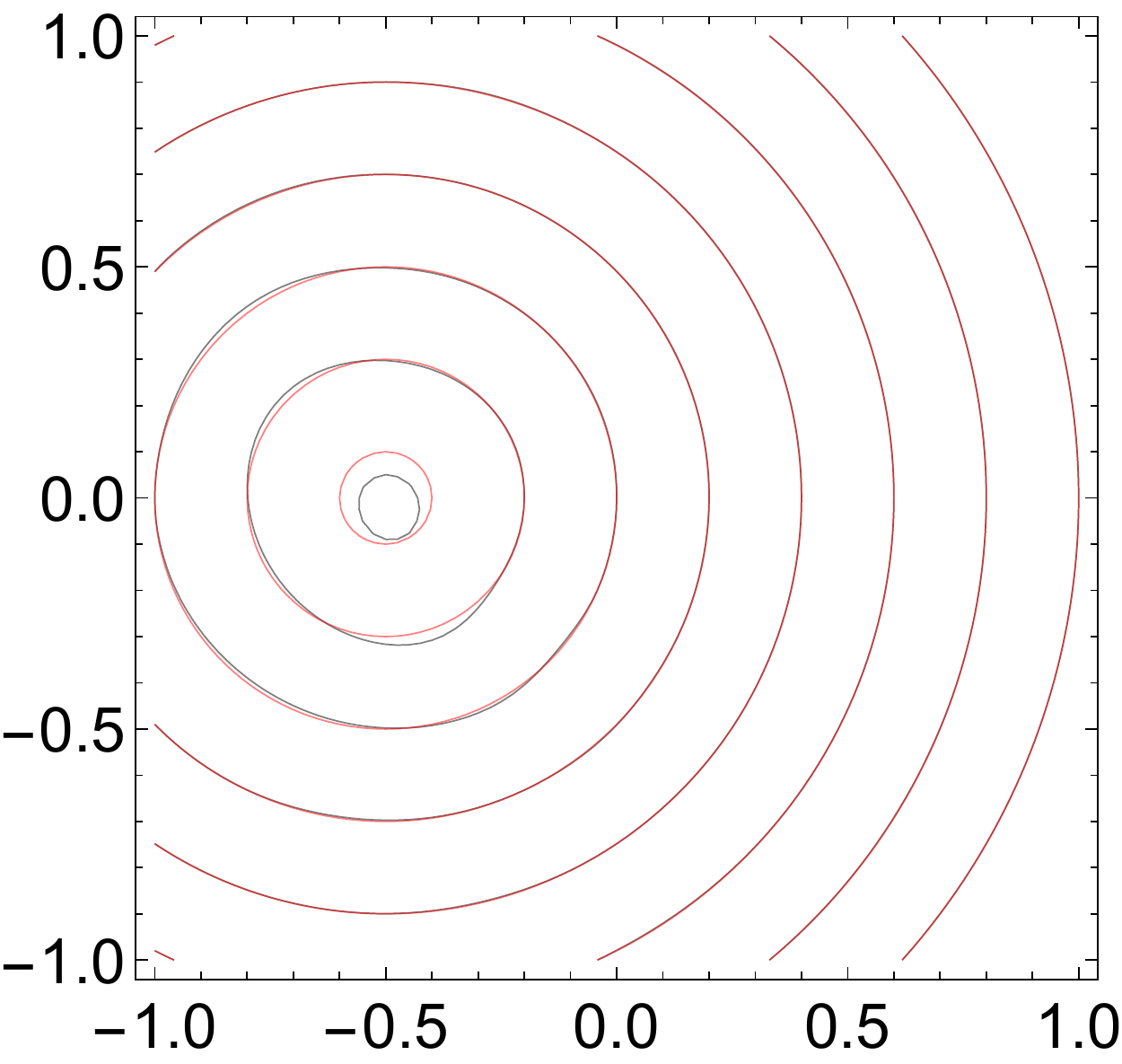}
\includegraphics[width=0.24\columnwidth]{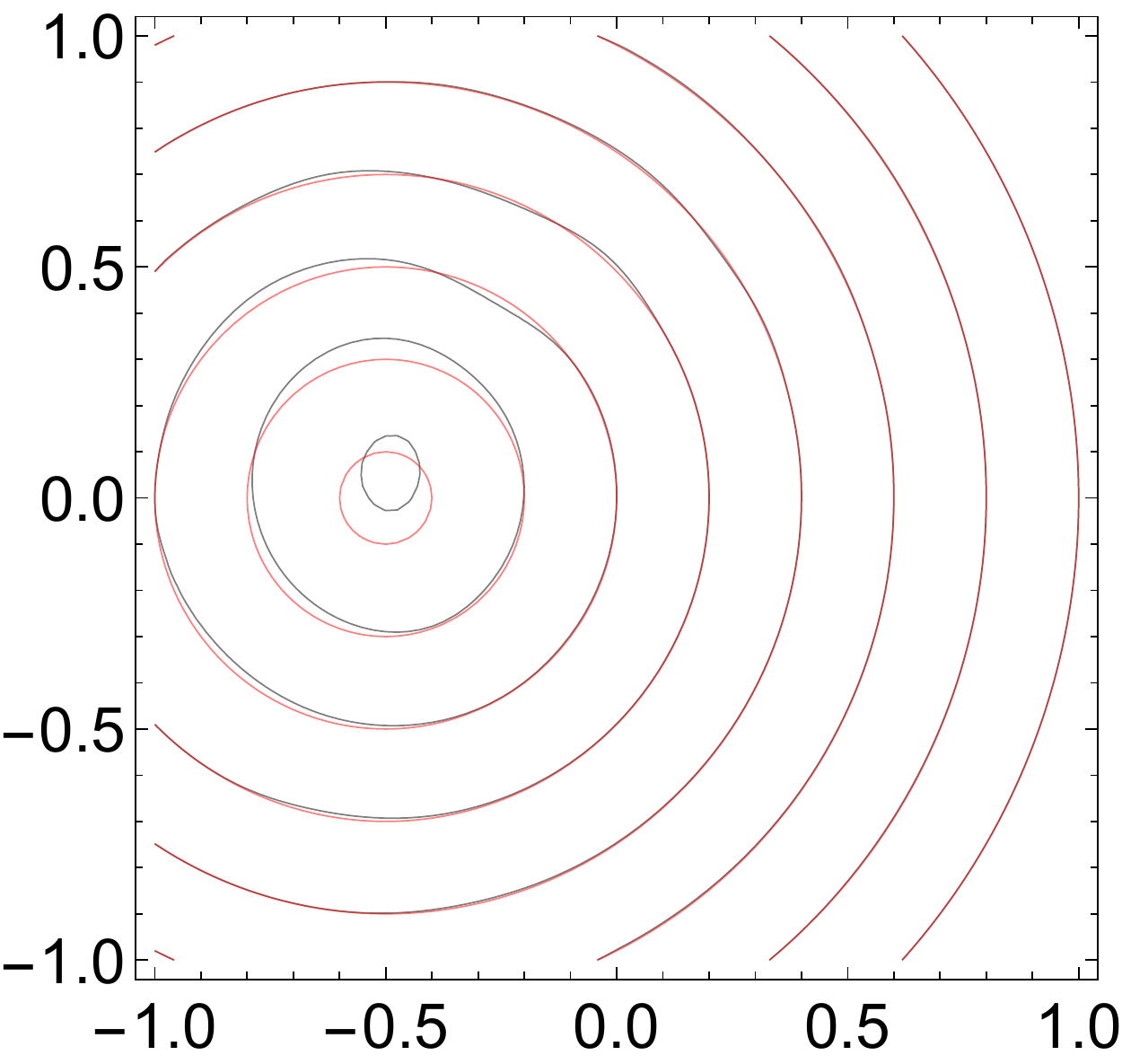}
\includegraphics[width=0.24\columnwidth]{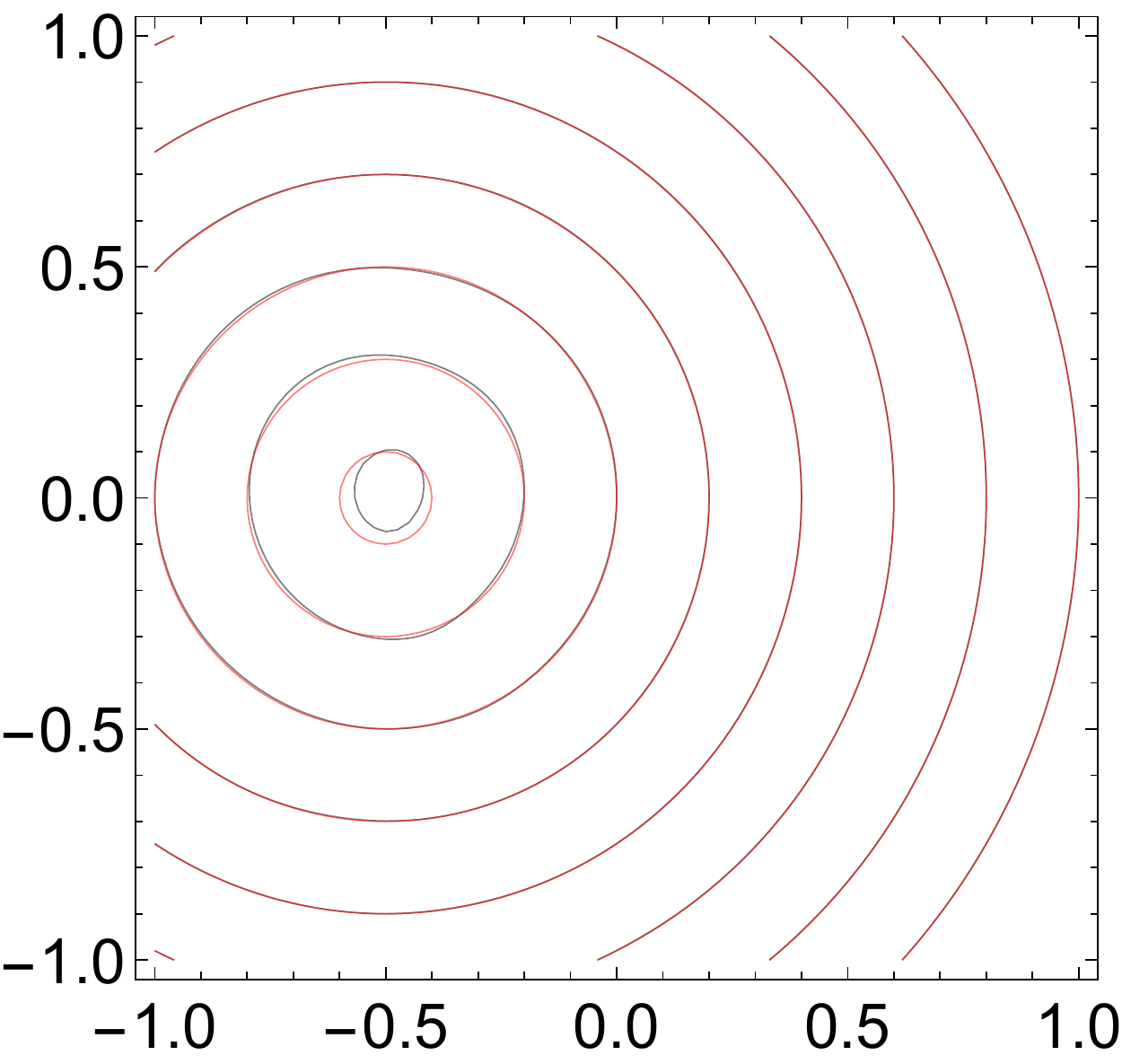}
\caption{Numerical solutions for the rotation of distance function for the coarse grid $M=40$ and small time step with $N=100$ (the first row of pictures) and for the medium grid $M=80$ and large time step with $N=40$ (the second row). The red contour lines represent the exact solution, the black contour lines represent the numerical solutions for values $0.1$, $0.3$ up to $1.5$ at $t=1$, see also the Table \ref{tab01}. The pictures represent the choices (\ref{kp}) - (\ref{k3}) in that order.
}
\label{fig01}
\end{center}
\end{figure}

In the next example we choose $u^0(x,y)= \sqrt{(x+0.5)^2+y^2}$, i.e. a distance function to point $(-0.5,0)$. The choices of $N$ and $M$ are identical to the previous example and the results are presented analogously to it  in Table \ref{tab01} and in the first row of pictures in Figure \ref{fig01}. One can observe the so-called ``phase error'' \cite{lev02} for the choices (\ref{kp}) and (\ref{km}) in Figure \ref{fig01} that is  reduced for the choice (\ref{k0}) in this example. 

\begin{table}[h!]
\begin{center}
\begin{tabular}{||c||c|c||c|c||c|c||c|c||c|c||}
\hline
$M$ & $E$ & $\min$  & $E$ & $\min$ & $E$ & $\min$ & $E$ & $\min$ & $E$ & $\min$ \\
\hline
40  & 32. & .10   & 26. & .13  & 4.4 & .091 & 3.8 & .087 & 4.0 & .087 \\
80  & 9.8 & .066 & 8.0 & .078 & 1.2 & .055 & .92 & .052 & .97 & .052\\
160& 2.9 & .042 & 2.4 & .049 & 3.2 & .033 & .23 & .031 & .23 & .031 \\
\hline
\end{tabular}\\[1ex]
\begin{tabular}{||c||c|c||c|c||c|c||c|c||c|c||}
\hline
$M$ & $E$ & $\min$  & $E$ & $\min$ & $E$ & $\min$ & $E$ & $\min$ & $E$ & $\min$ \\
\hline
40  & 124. & .16   & 23.3 & .13  & 61.7 & .13 & 21.8 & .12 & 20.1 & 0.12 \\
80  & 41.3 & .11 & 8.17 & .078 & 19.7 & .082 & 6.90 & .075 & 6.10 & .072 \\
160& 12.6 & .066 & 2.52 & .049 & 5.90 & .051 & 2.05 & .046 & 1.67 & .043 \\
\hline
\end{tabular}
\caption{The error $E$ in (\ref{error}) (multiplied by $10^3$) and the minimum of numerical solutions for the rotation of distance function using small time steps with $N=5 M/2$ (the top table) and large time steps with  $N=M/2$ (the bottom table). The results are obtained
using (\ref{si2d}) with (\ref{kp})  (the $2^{nd}$ - $3^{rd}$ columns), (\ref{km}) (the $4^{th}$ - $5^{th}$ ones), (\ref{k0}) (the $6^{th}$ - $7^{th}$ ones), (\ref{k3}) (the $8^{th}$ - $9^{th}$ ones), and (\ref{ctusi2d1}) with (\ref{k3}). }
\label{tab01}
\end{center}
\end{table}

To illustrate also the numerical stability conditions of schemes we recompute the example with $5$ times larger time step, i.e. $N=M/2$ when the maximum of Courant numbers $|\mathcal{C}_{ij}|$ and $|\mathcal{D}_{ij}|$ is approximately $6.28$. The results are given again in Table \ref{tab01} and in the second row of pictures in Figure \ref{fig01}. All schemes give stable numerical solutions 
with the CTU extension bringing the lowest error $E$ for this setting. One can note an interesting fact that the variant (\ref{km}) gives  for larger time step a lower error $E$ than the variant (\ref{k0}) in this example.

\begin{figure}[h!]
\begin{center}
\includegraphics[width=0.24\columnwidth]{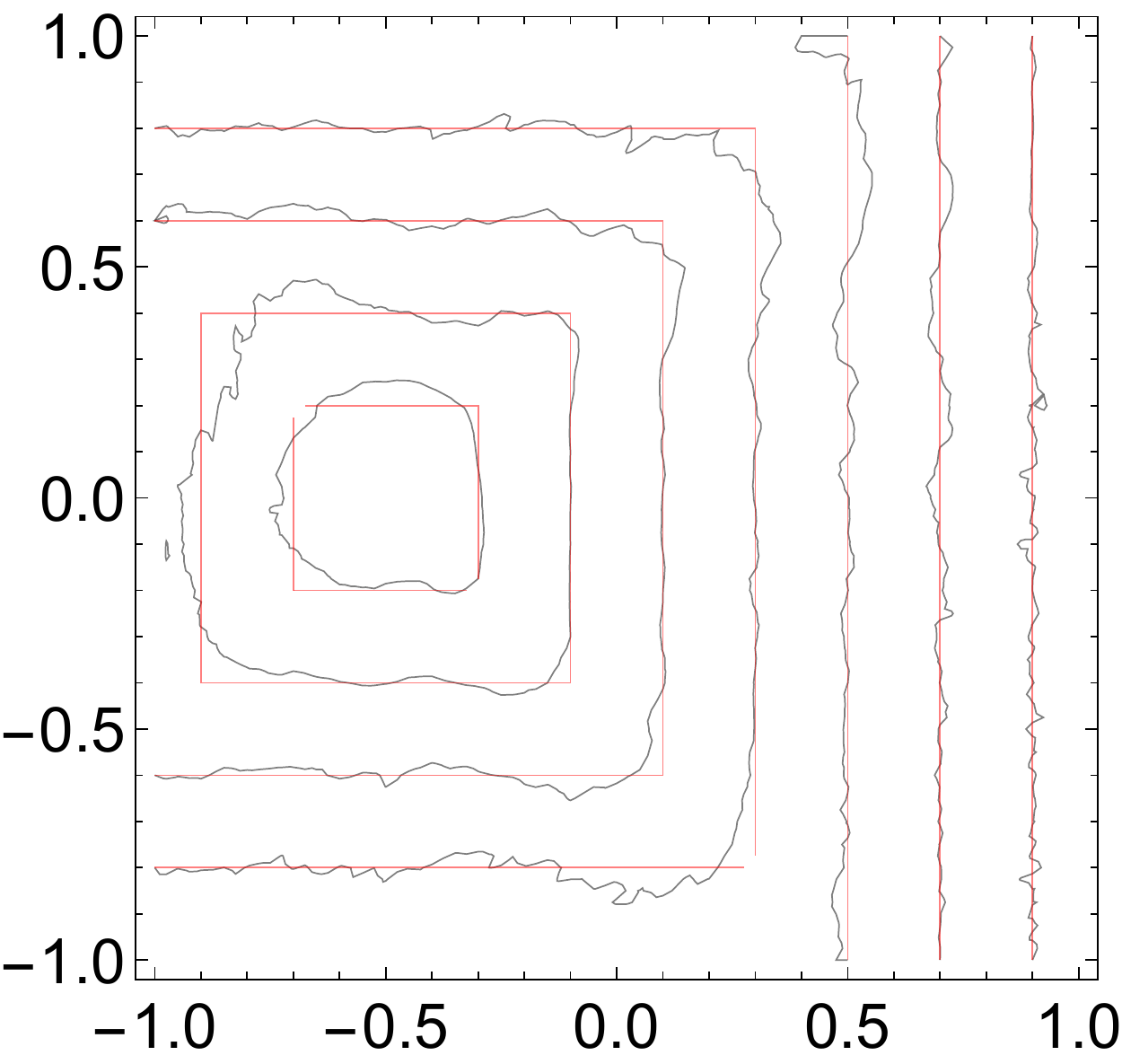}
\includegraphics[width=0.24\columnwidth]{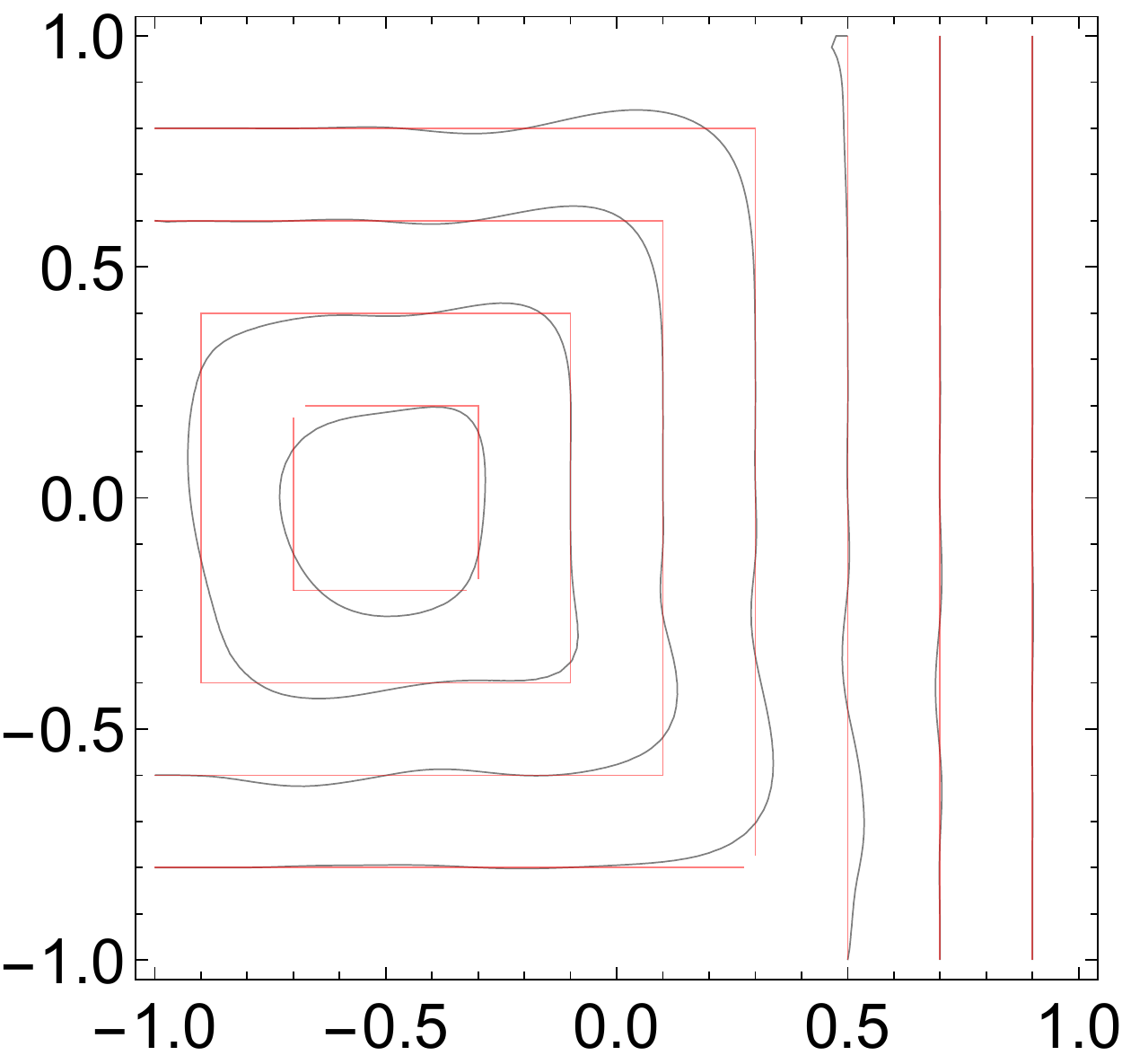}
\includegraphics[width=0.24\columnwidth]{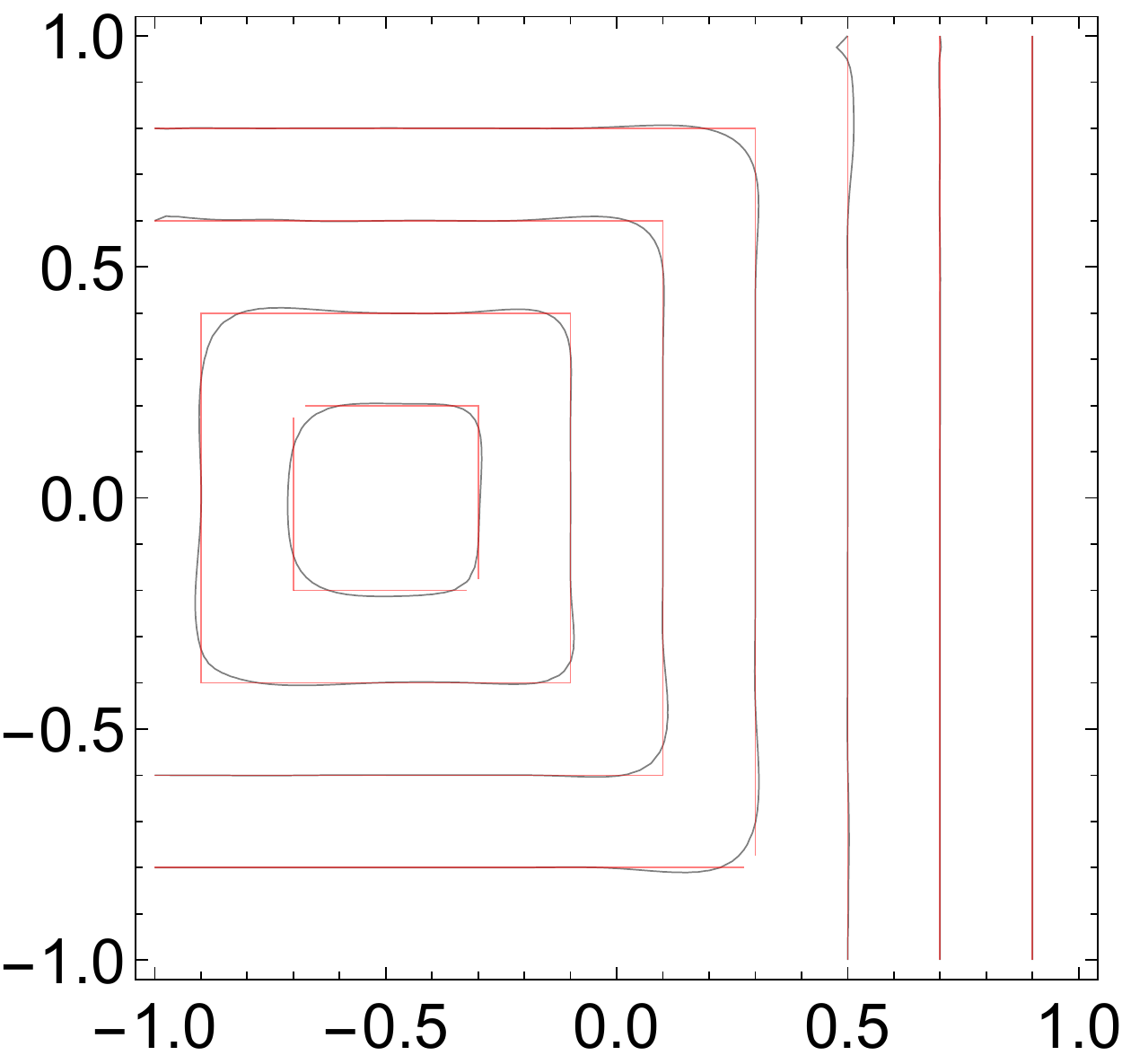}
\includegraphics[width=0.24\columnwidth]{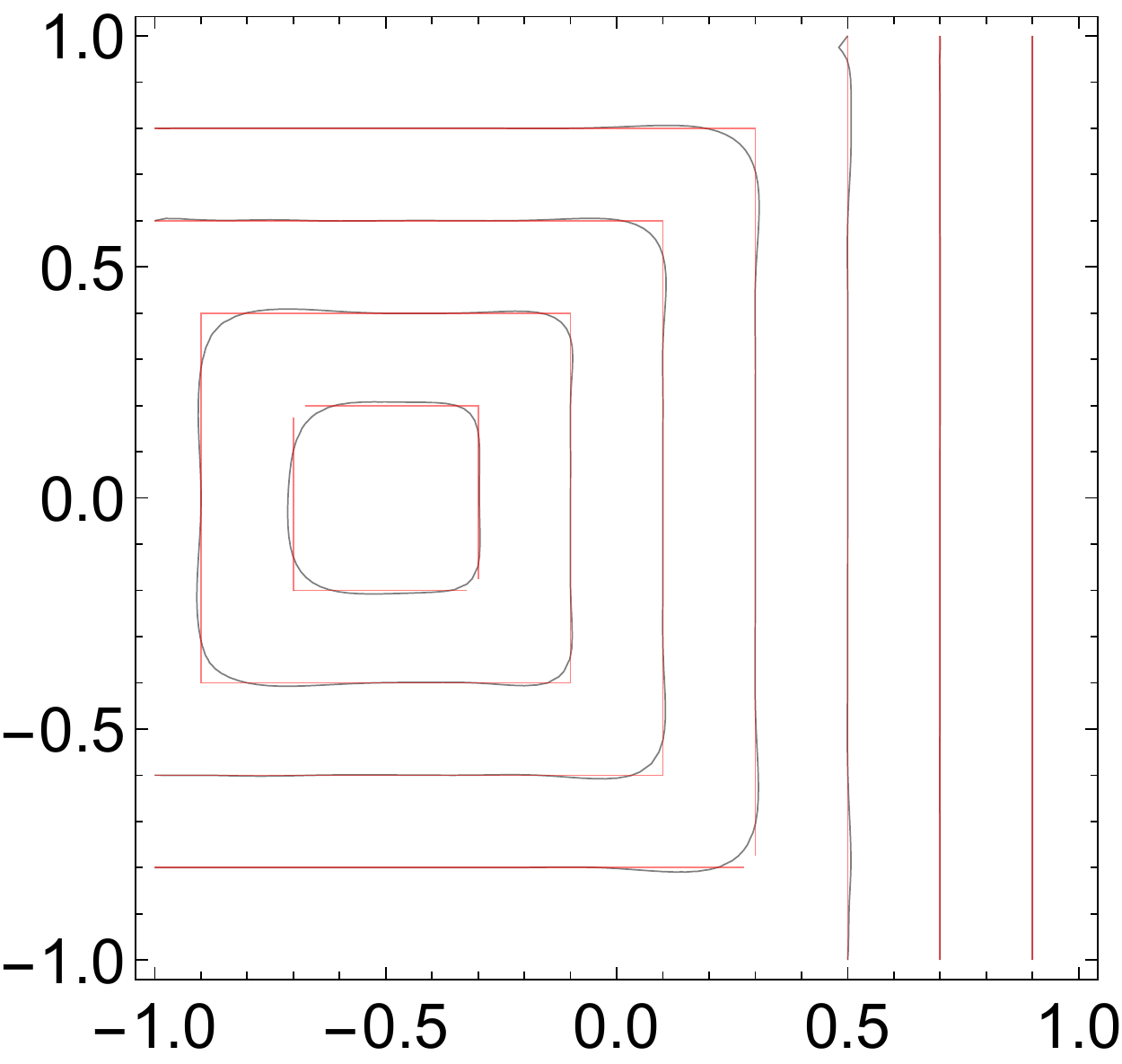}
\caption{Numerical solutions for the rotation of maximum distance function for the medium grid $M=80$ and small time step with $N=100$. The red contour lines represent the exact solution, the black contour lines represent the numerical solutions for values $0.2$, $0.4$ up to $1.4$ at $t=1$. The pictures represent the scheme (\ref{si2d}) with the choices (\ref{kp}) - (\ref{k0}) and (\ref{ctusi2d1}) with (\ref{k3}) in that order.
}
\label{fig02}
\end{center}
\end{figure}

Next we choose $u^0(x,y)= \max\{(x+0.5),y\}$, i.e. a distance function in maximum metric to point $(-0.5,0)$. 
The results are summarized in Table \ref{tab02} and in Figure \ref{fig02}.  All schemes seem to show the EOC having better than the $1^{st}$ order accuracy. The scheme (\ref{si2d}) with the choice (\ref{kp}), for which the amplification factor is constant and equals $1$, produces oscillations that are not damped by the scheme.

\begin{table}[h!]
\begin{center}
\begin{tabular}{||c||c|c||c|c||c|c||c|c||c|c||}
\hline
$M$ & $E$ & $\min$  & $E$ & $\min$ & $E$ & $\min$ & $E$ & $\min$ & $E$ & $\min$ \\
\hline
40  & 115.& 9.2 & 75. & 11.& 28. & 8.2 & 28. & 7.9 & 26. & 7.9 \\
80  & 62. & 6.4 & 38. & 7.3 & 11. & 5.0 & 11. & 4.8 & 10. & 4.8\\
160& 30. & 3.5 & 17. & 46. & .44 & 3.0 & .42 & 2.8 & .39 & 2.8 \\
\hline
\end{tabular}
\caption{The error (\ref{error}) (multiplied by $10^3$) and the minimum of numerical solutions for the rotation of distance function using (\ref{kp})  (the $2^{nd}$ - $3^{rd}$ columns), (\ref{km}) (the $4^{th}$ - $5^{th}$ ones), (\ref{k0}) (the $6^{th}$ - $7^{th}$ ones), and (\ref{k3}).  Note that $N=5/2 M$.}
\label{tab02}
\end{center}
\end{table}

Finally, we choose the single vortex example that is characterized by strongly deformable flow. The velocity $\vec{V}=(V,W)$ is given by
\begin{eqnarray*}
V(x,y)=-4 \sin^2(\pi(x+1)/2) \sin(\pi(y+1)/2) \cos(\pi(y+1)/2) \,,\\
W(x,y)=4 \sin^2(\pi(y+1)/2) \sin(\pi(x+1)/2) \cos(\pi(x+1)/2) \,.
\end{eqnarray*}
for $t \in [0,2.5]$. The initial function $u^0(x,y)$ is a distance function to a circle of radius $0.3$ with the center in $(0,0.5)$. As $\vec{V} \equiv 0$ at the boundary $\partial D$ of square domain $D=(-1,1)^2$, one can fix the values of $u(x,y,t)=u^0(x,y)$ for $(x,y) \in \partial D$ and $t \in [0,2.5]$.

We present the results at $t=2.5$ where the largest deformation of initial function can be observed. As the exact solution at this time is not available, we compute a reference numerical solution $\tilde U_{i j}^N$ obtained with  $N=5 {\cal M}/4$ and ${\cal M}=1280$ and compare it with numerical solutions for $M=80,160,320$ at $N=5 M/4$ using 
\begin{equation}
\label{l1e}
e = \frac{4}{M^2} \sum_{i,j=0}^M  |U_{i j }^{5M/4} - \tilde U_{{\cal M}/M i \, {\cal M}/M j }^{5{\cal M}/4}| \,.
\end{equation}
The chosen time step corresponds to the maximum of $|\mathcal{C}_{ij}|$ and $|\mathcal{D}_{ij}|$ being 2.

\begin{figure}[h!]
\begin{center}
\includegraphics[width=0.24\columnwidth]{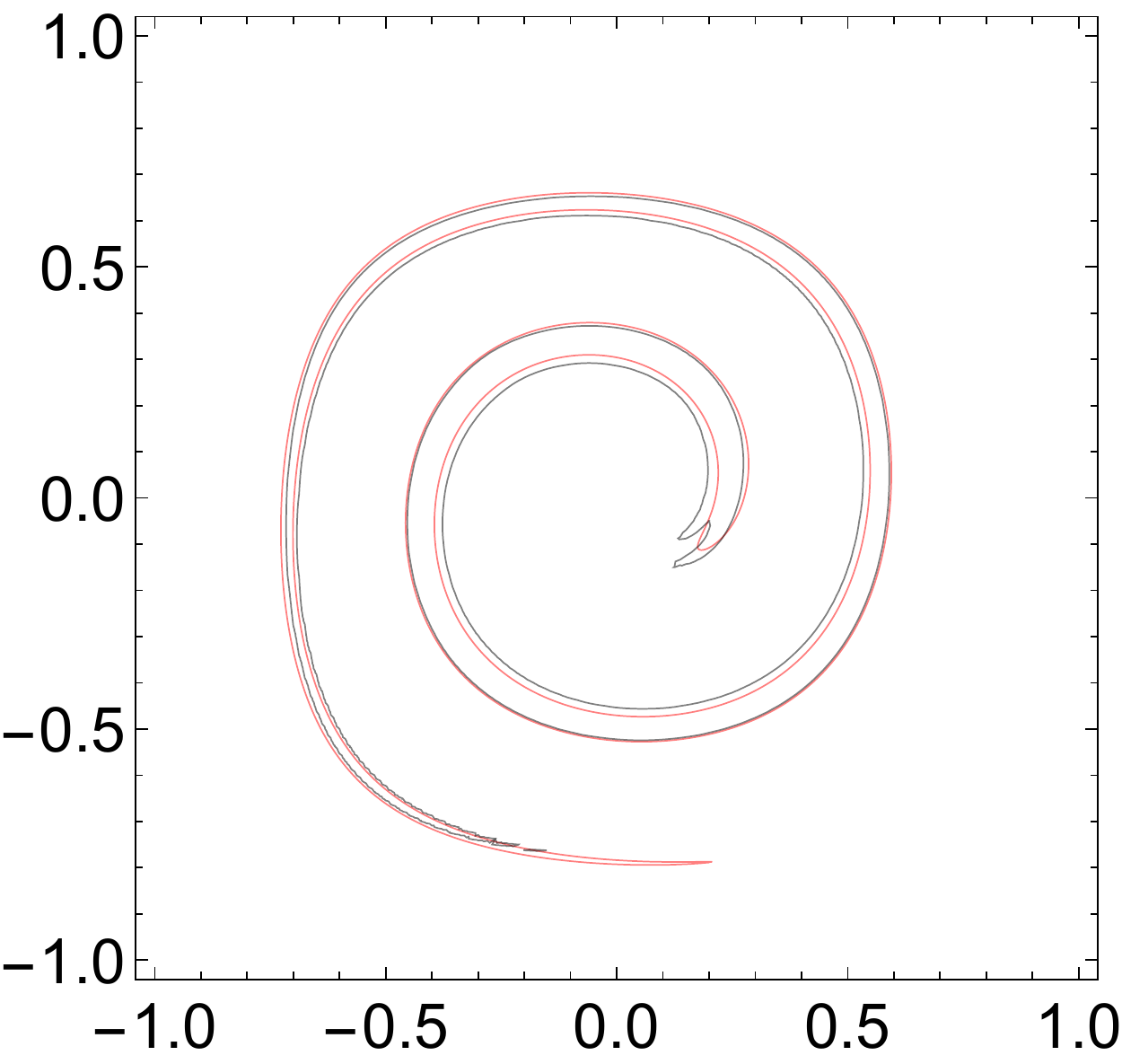}
\includegraphics[width=0.24\columnwidth]{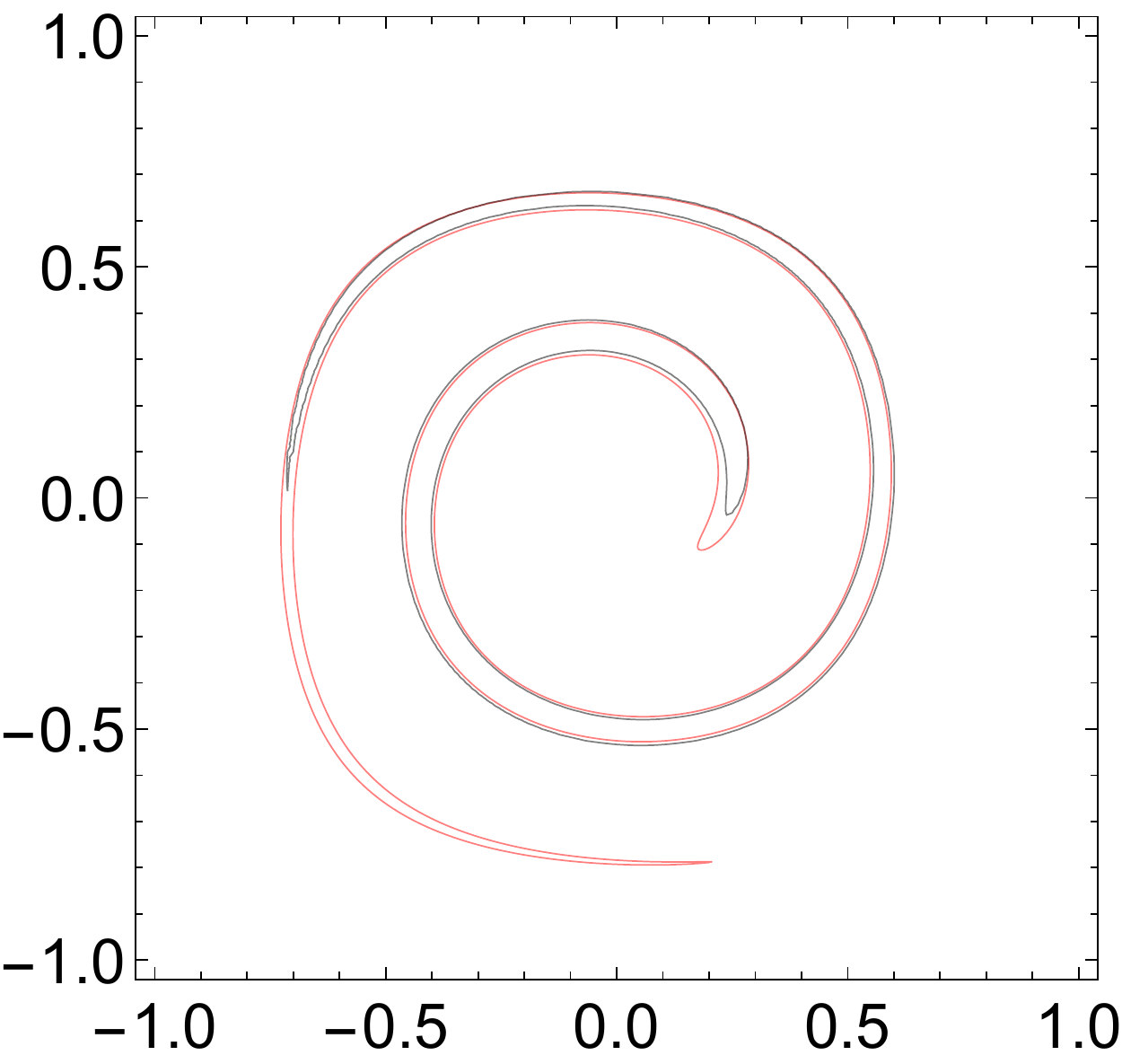}
\includegraphics[width=0.24\columnwidth]{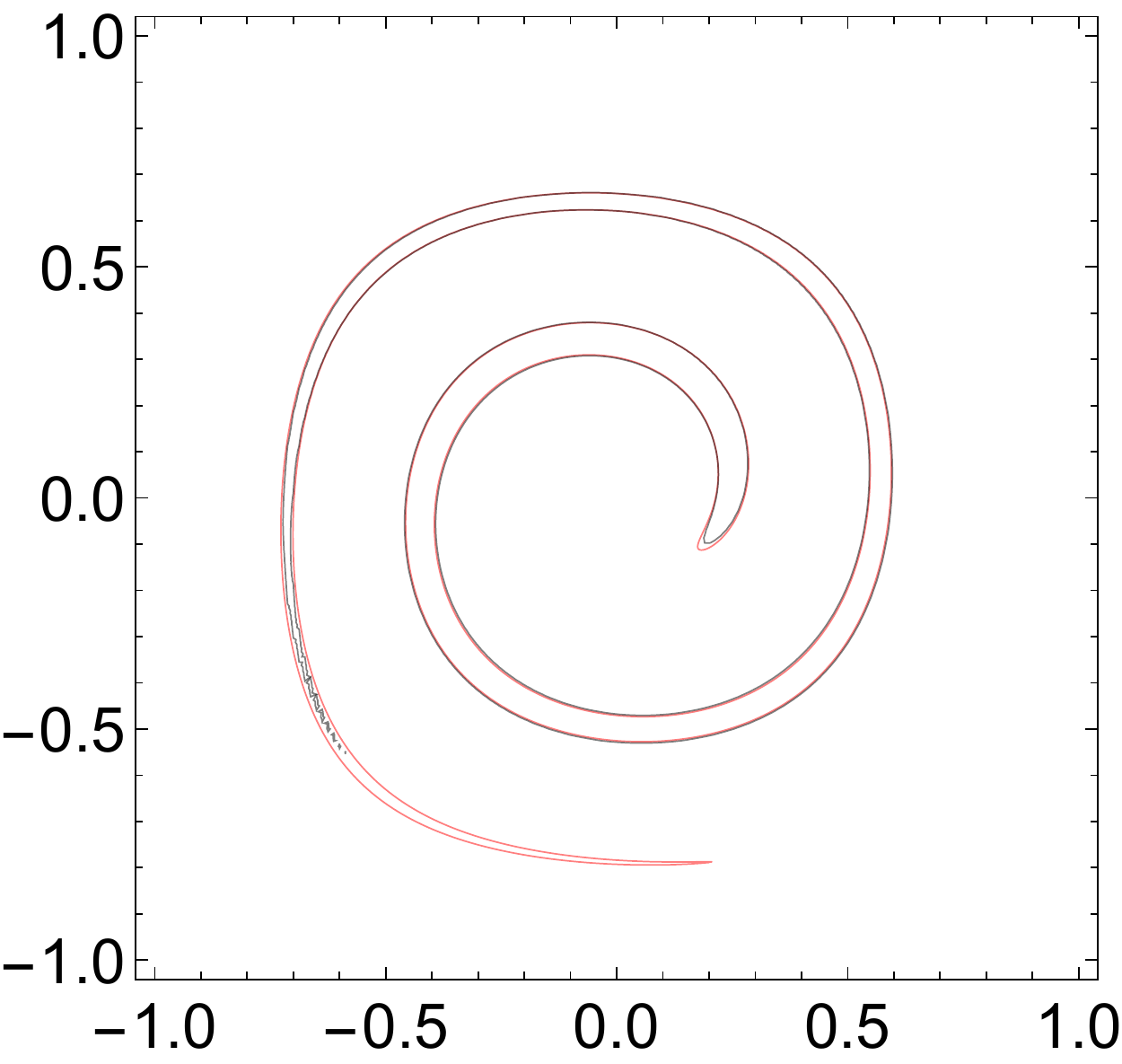}
\includegraphics[width=0.24\columnwidth]{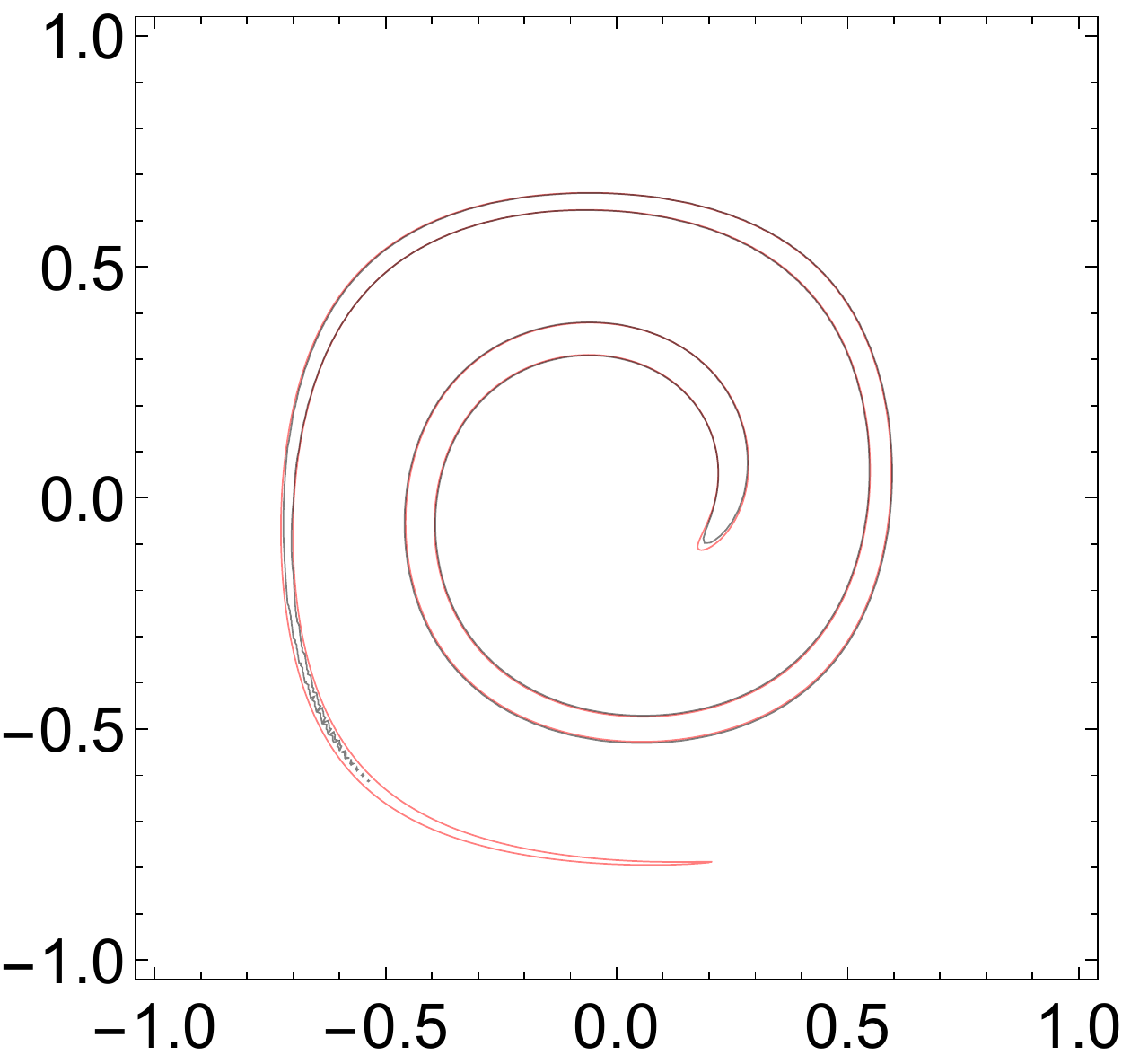}
\caption{Numerical solutions at $t=2.5$ for the single vortex example for the medium grid $M=160$ and $N=5 M/4$. The red contour line in each plot represents the contour line for the value $0.3$ of reference solution, the black contour lines represent it for the numerical solutions of scheme (\ref{si2d}) with the variants (\ref{kp}) - (\ref{k3}) in that order.
}
\label{fig03}
\end{center}
\end{figure}

We use the Dirichlet boundary conditions and we set $\kappa^x=1$ or $\kappa^y=1$ for points $(x_i,y_j)$ at the neighborhood of boundary nodes, i.e. $i=1$ or $i=M-1$ for arbitrary $j$ and so on, to avoid in the implicit part of scheme (\ref{si2d}) the values of numerical solution in points lying outside of the domain. The results are summarized in Table \ref{tab03} and the visual comparison is given in Figure \ref{fig03}. 

\begin{table}[h!]
\begin{center}
\begin{tabular}{||c||c|c||c|c||c|c||c|c||c|c||}
\hline
$M$ & $e$ & $\min$  & $e$ & $\min$ & $e$ & $\min$ & $e$ & $\min$ & $e$ & $\min$ \\
\hline
80  & 1.86 & 1.10 & 1.60 & 3.03 & .814 & 2.16 & .800 & 2.25 & .894 & 2.36 \\
160 & .743& 1.21 & .650 & 1.81 & .301 & 1.34 & .282 & 1.38 & .301 & 1.44 \\
320& .256 & .794 & .254 & 1.12 & .085 & .815 & .077 & .830 & .082 & .851 \\
\hline
\end{tabular}
\caption{The error (\ref{l1e}) (multiplied by $10$) and the minimal value (multiplied by $10$) for numerical solutions of single vortex example that is solved using (\ref{si2d}) with (\ref{kp})  (the $2^{nd}$ - $3^{rd}$ columns), (\ref{km}) (the $4^{th}$ - $5^{th}$ ones), (\ref{k0}) (the $6^{th}$ - $7^{th}$ ones),  (\ref{k3}) (the $8^{th}$ - $9^{th}$ ones)  and with (\ref{ctusi2d1}) using (\ref{k3}). Note that $N=5 M/4$.}
\label{tab03}
\end{center}
\end{table}

\begin{figure}[h!]
\begin{center}
\includegraphics[width=0.3\columnwidth]{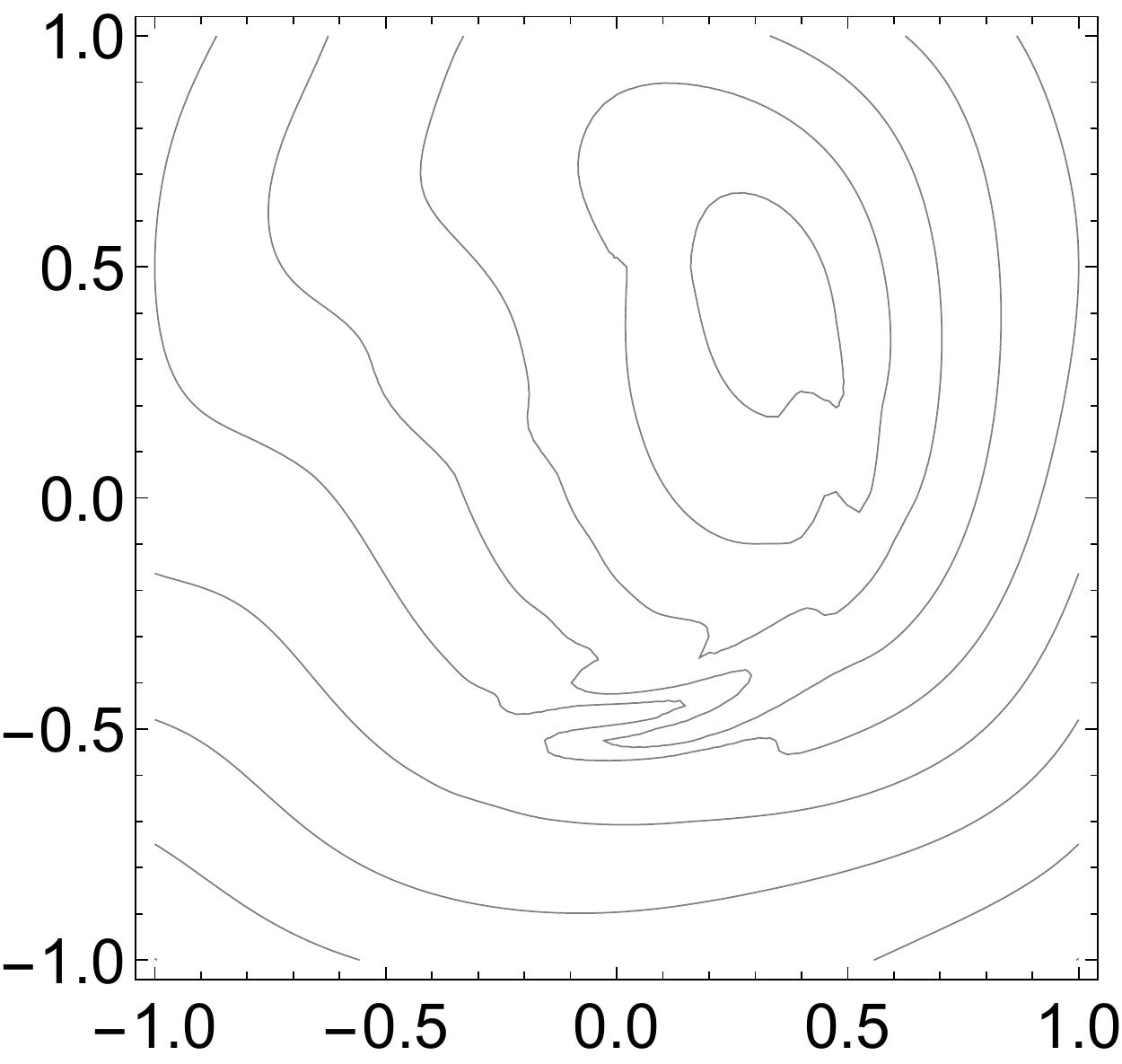}
\includegraphics[width=0.3\columnwidth]{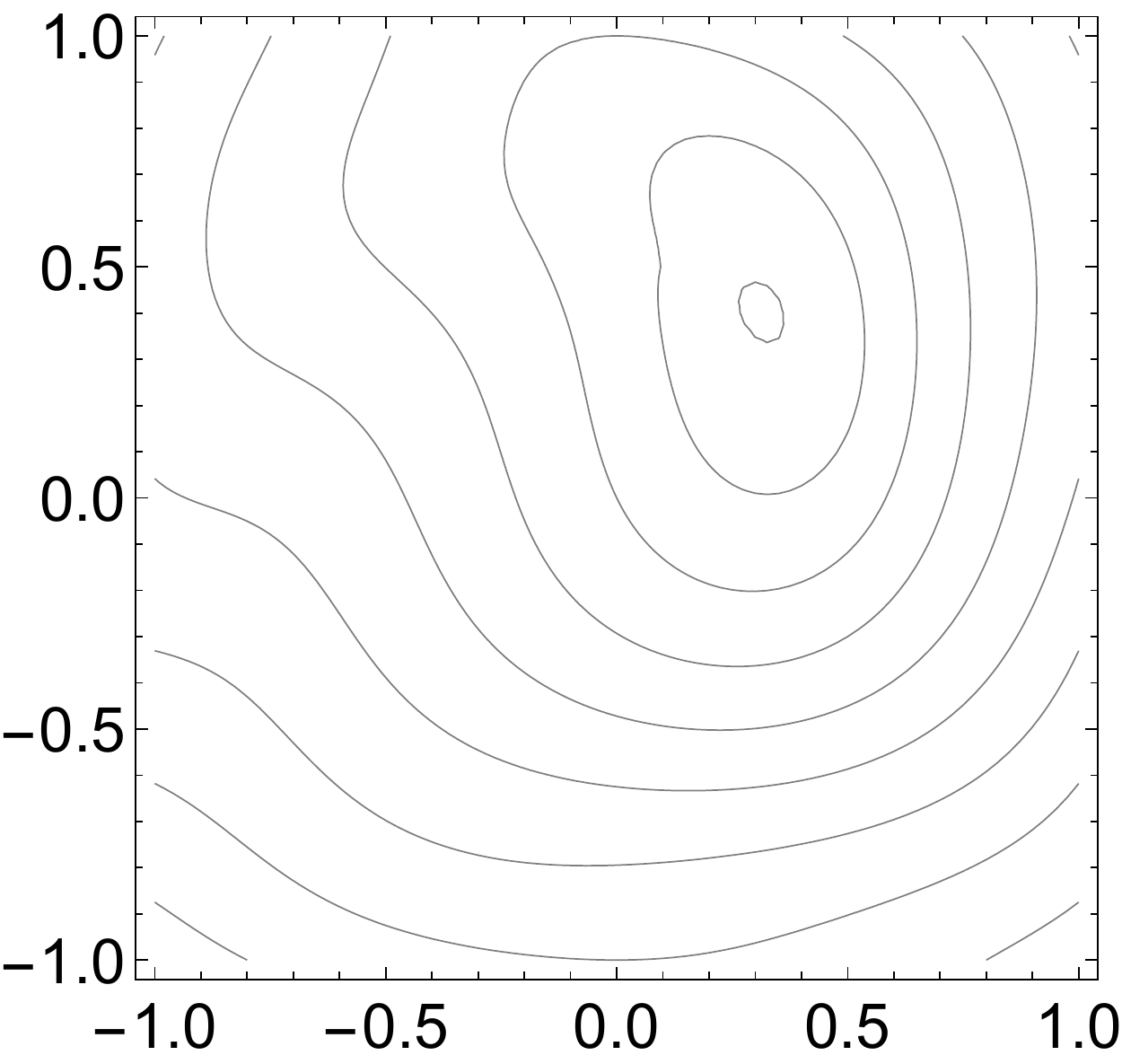}
\includegraphics[width=0.3\columnwidth]{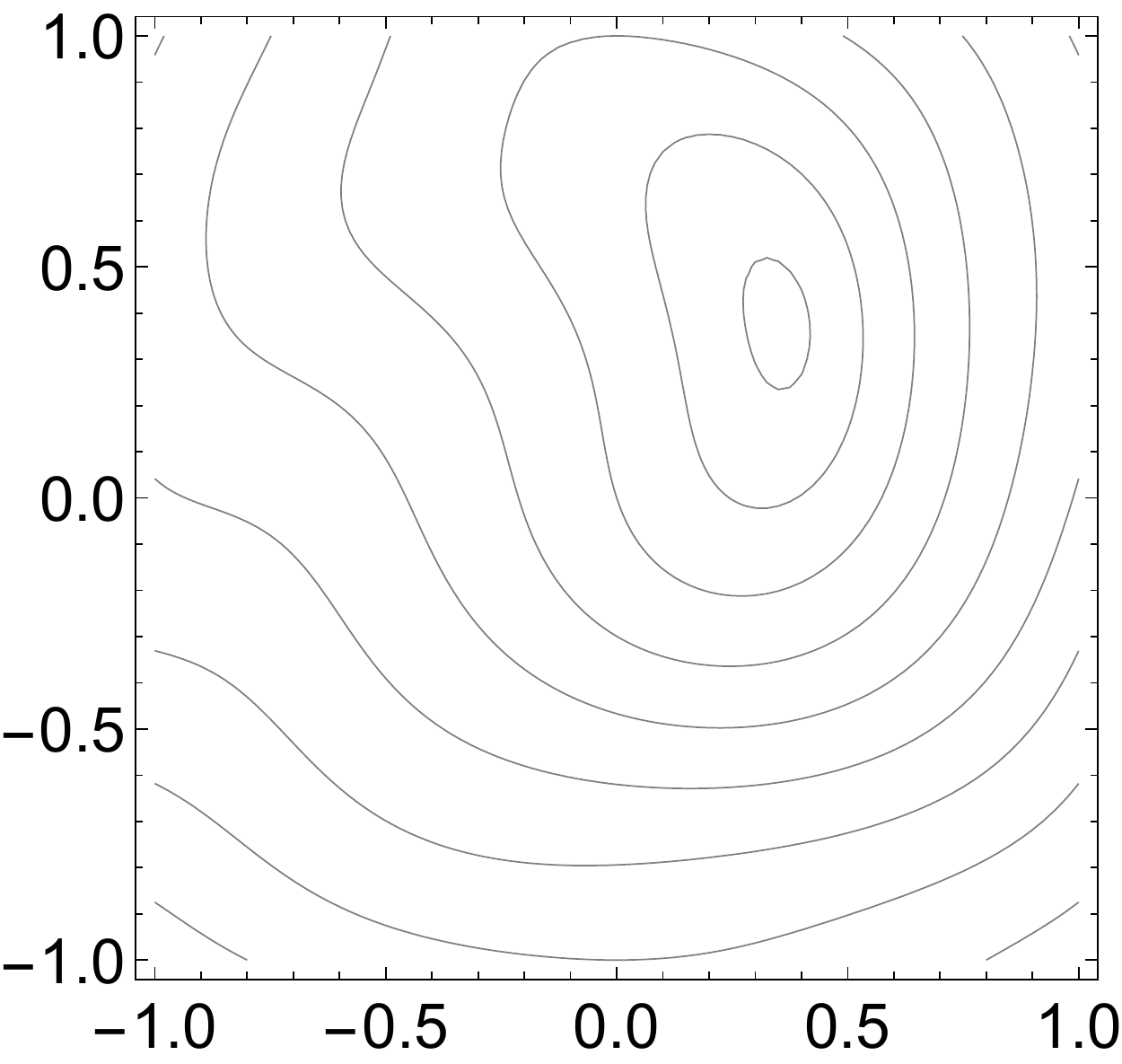}
\caption{Numerical solutions at $t=0.2$ with contour lines for values $0.1$, $0.3$ up to $1.7$ for the single vortex example for the coarse grid $M=80$. The first picture uses one time step (the Courant number equals $16$) with the scheme (\ref{si2d}) using (\ref{k3}) when instabilities can be observed. The second picture is obtained using the CTU extension (\ref{ctusi2d1}) with (\ref{k3}) taking one identical time step when no instabilities occur. For a comparison the last picture is obtained using $16$ uniform time steps that corresponds the maximal Courant number being $1$ using (\ref{si2d}) with (\ref{k3}).
}
\label{fig04}
\end{center}
\end{figure}

Concerning numerical stability analysis we can confirm our results discussed in section \ref{sec-2d} for this example. We compute one large time step with the scheme (\ref{si2d}) and (\ref{k3})  such that the maximum of Courant numbers $|\mathcal{C}_{ij}|$ and $|\mathcal{D}_{ij}|$  equals $16$ for $M=80$. In Figure \ref{fig04} we present the numerical result obtained after one sweep in iterative solver (i.e. four Gauss-Seidel iterations) where instabilities can be clearly observed and that grow when more sweeps are used. Computing the example with identical one time step, but using the CTU extension (\ref{ctusi2d1}) with (\ref{k3}), no instabilities occur. We note that instabilities can be observed also when using (\ref{si2d}) with (\ref{km}) although even larger time step and larger number of sweeps must be used.

\section{Conclusions}
\label{sec-con}

In this paper a class of semi-implicit schemes on Cartesian grids for the numerical solution of linear advection equation is derived. The derivation follows the partial Cauchy-Kowalewski (or Lax-Wendroff) procedure to replace the time derivatives of exact solution in Taylor series. Opposite to the standard form of this procedure when only space derivatives are used in the replacement, the partial procedure exploits also the mixed time-space derivatives. 

Firstly, we derive the one dimensional semi-implicit $\kappa$-scheme that is $2^{nd}$ order accurate with unconditional numerical stability for variable velocity and for all considered values of $\kappa$. The analogous fully implicit $\kappa$-scheme has unconditional numerical stability only for a limited range of $\kappa$ and it can be shown only for constant velocity case. 

We note that the fully explicit, fully implicit and semi-implicit $\kappa$-schemes offer a value of $\kappa$ (in each case different one) for which the scheme is $3^{rd}$ order accurate in time and space for constant velocity. Nevertheless only the semi-implicit variant has unconditional numerical stability for this value.

The dimension by dimension extension of one dimensional semi-implicit $\kappa$-scheme to finite difference grids in several dimensions gives a $2^{nd}$ order accurate scheme for variable velocity case. This is a unique property of the semi-implicit scheme as analogous extensions of fully explicit or fully implicit $\kappa$-scheme result in only first order accurate schemes. 

The numerical von Neumann stability analysis for general settings of such semi-implicit $\kappa$-scheme gives only conditional stability. Moreover, the dimension by dimension extension of semi-implicit $\kappa$-scheme does not offer a choice of $\kappa$ parameters for which the scheme would be $3^{rd}$ order accurate for constant velocity case.

To overcome these disadvantages we derive the Corner Transport Upwind extension of two-dimensional semi-implicit $\kappa$-scheme by extending its stencil using diagonal upwind corner values in the implicit part. The resulting scheme has all desired properties that are considered in this paper - the scheme is $2^{nd}$ order accurate with numerically suggested unconditional stability for variable velocity case and it is $3^{rd}$ order accurate if the velocity is constant. 

 
\bibliographystyle{plain}
\bibliography{lit}

\end{document}